\documentclass[MSNbibl,number,citesort,seceqn,dvips]{arxbj}
\usepackage{mathbh,upgreek}

\aid{0}
\volume{18}
\issue{2}
\pubyear{2012}
\firstpage{552}
\lastpage{585}
\doi{10.3150/11-BEJ354}

\makeatletter
\newcommand{\iint}{\int\!\!\!\int}
\newtheorem{thmm}{Theorem}[section]
\newtheorem{lem}{Lemma}[section]
\newproclaim{defn}{Definition}[section]

\newremark{rem}{Remark}

\def\N{\mathbb{N}}
\def\R{\mathbb{R}}
\def\C{\mathbb{C}}
\def\Z{\mathbb{Z}}
\def\E{\mathbb{E}}
\def\I{\mathbh{1}}
\def\X{\mathbb{X}}
\def\XX{\mathfrak{X}}
\def\ninfty{\mathop{\longrightarrow}_{n\to\infty}}
\makeatother

\begin{document}
\begin{frontmatter}

\title{Degenerate $U$- and $V$-statistics under weak dependence:
Asymptotic theory and bootstrap consistency}
\runtitle{Bootstrap for $U$-statistics under weak dependence}

\begin{aug}
\author{\fnms{Anne} \snm{Leucht}\corref{}\ead[label=e1]{anne.leucht@uni-jena.de}}
\runauthor{A. Leucht}
\address{Friedrich-Schiller-Universit\"at Jena,
Institut f\"ur Stochastik,
Ernst-Abbe-Platz 2,
D-07743 Jena,
Germany. \printead{e1}}
\end{aug}

\received{\smonth{2} \syear{2010}}
\revised{\smonth{11} \syear{2010}}

%
\begin{abstract}
We devise a general result on the consistency of model-based bootstrap
methods for $U$- and $V$-statistics under easily verifiable conditions.
For that purpose, we derive the limit distributions of degree-2
degenerate $U$- and $V$-statistics for weakly dependent $\R^d$-valued
random variables first. To this end, only some moment conditions and
smoothness assumptions concerning the kernel are required. Based on
this result, we verify that the bootstrap counterparts of these
statistics have the same limit distributions. Finally, some
applications to hypothesis testing are presented.
\end{abstract}

%

\end{frontmatter}

\section{Introduction}
\label{S1}

Numerous test statistics can be formulated or approximated in terms of
degenerate $U$- or $V$-type statistics. Examples include the Cram\'er--von
Mises statistic, the Anderson--Darling statistic or the $\chi
^2$-statistic.
For i.i.d.~random variables the limit distributions of $U$- and
$V$-statistics can be derived via a spectral decomposition of their
kernel if the latter is squared integrable.
To use the same method for dependent data, often restrictive
assumptions are required whose validity is quite complicated or even
impossible to verify in many cases. The first of our two main results
is the derivation of the asymptotic distributions of $U$- and
$V$-statistics under assumptions that are fairly easy to check. This
approach is based on a wavelet decomposition instead of a spectral
decomposition of the kernel.

The limit distributions for both independent and dependent observations
depend on certain parameters which in turn depend on the underlying
situation in a complicated way. Therefore, problems arise as soon as
critical values for test statistics of $U$- and $V$-type have to be
determined. The bootstrap offers a convenient way to circumvent these
problems; see Arcones and Gin{\'e} \cite{ArGi92}, Dehling and Mikosch \cite{DM94} or
Leucht and Neumann~\cite{LeuNeu08} for the
i.i.d.~case. To our knowledge, there are no results concerning
bootstrapping general degenerate $U$-statistics of non-independent
observations. As a second main result of the paper, we establish
consistency of model-based bootstrap methods for $U$- and $V$-type
statistics of weakly dependent data.

In order to describe the dependence structure of the sample, we do not
invoke the concept of mixing although a great variety of processes
satisfy these constraints and various tools of probability theory and
statistics such as central limit theorems, probability and moment
inequalities can be carried over from the i.i.d.~setting to mixing processes.
However, these methods of measuring dependencies are inappropriate in
the present context since not only the asymptotic
behaviour of $U$- and $V$-type statistics but also bootstrap
consistency is focused. Model-based bootstrap methods can yield samples
that are no longer mixing even though the original sample satisfies
some mixing condition. A simple example is presented in Section~\ref
{SS42}. There we consider a model-specification test within the class
of nonlinear $\operatorname{AR}(1)$ processes. Under ${\mathcal H_0}$,
$X_k=g_0(X_{k-1})+\varepsilon_k$, where $g_0$ is Lipschitz contracting
and $(\varepsilon_k)_k$ is a sequence of i.i.d.~centered innovations.
It is most natural to draw the bootstrap innovations $(\varepsilon
_k^*)_k$ via Efron's bootstrap from the
recentered residuals first. Then the bootstrap counterpart of $(X_k)_k$
is generated iteratively by choosing an initial variable~$X_0^*$
independently of $(\varepsilon_k^*)_k$ and defining
$X_k^*=g_0(X_{k-1}^*)+\varepsilon_k^*$. Due to the discreteness of the
bootstrap innovations, commonly used coupling techniques to prove
mixing properties for Markovian processes fail; see also Andrews \cite{An84}.
It turns out that the characterization of dependence structures
introduced by Dedecker and Prieur~\cite{DP05} is exceptionally suitable here. Based on
their $\tau$-dependence coefficient it is possible to construct an
$L_1$-coupling in the following sense. Let ${\mathcal M}$ denote a
$\sigma$-algebra generated by sample variables of the ``past'' and let
$X$ be a random variable of a certain ``future'' time point. Then, the
minimal $L_1$-distance between $X$ and a random variable that has the
same distribution as $X$ but that is independent of ${\mathcal M}$ is
equivalent to the $\tau$-dependence coefficient $\tau({\mathcal M}, X)$.

We exploit this coupling property in order to derive the asymptotic
distribution for the original as well as the bootstrap statistics of
degenerate $U$-type. Basically, both proofs follow the same lines.
First, the (almost) Lipschitz continuous kernels of the $U$-statistics
are approximated by a finite wavelet series expansion. There are two
crucial points that assure asymptotic negligibility of the
approximation error. On the one hand, the smoothness of the kernel
function carries over to its wavelet approximation uniformly in scale,
cf.~Lemma~\ref{l.2}. On the other hand, Lipschitz continuity of the
kernel and the $L_1$-coupling property of the underlying $\tau
$-dependent sample perfectly fit together. A next step contains the
application of a central limit theorem and the continuous mapping
theorem to determine the limits of the approximating statistics of
$U$-type. Based on these investigations, the asymptotic distribution of
the $U$-statistic and its bootstrap counterpart is then deduced via
passage to the limit. It can be expressed as an infinite weighted sum
of normal variables.

Our paper is organized as follows.
We start with an overview of asymptotic results on degenerate $U$-type
statistics of dependent random variables. In Section~\ref{SS22}, we
introduce the underlying concept of weak dependence and derive the
asymptotic distributions of $U$- and $V$-statistics.
On the basis of these results, we deduce consistency of general
bootstrap methods in Section~\ref{S3}. Some applications of the theory
to hypothesis testing are presented in Section~\ref{S4}.
All proofs are deferred to a final Section~\ref{S5}.

\section{Asymptotic distributions of $U$- and $V$-statistics}\label{S2}

\subsection{Survey of literature}\label{SS21}

Let $(X_n)_{n\in\N}$ be a sequence of $\R^d$-valued random variables
with common distribution~$P_X$. In the case of i.i.d.~random variables,
the limit distributions of degenerate $U$- and $V$-type statistics,
that is,
\[
n U_n=\frac{1}{n}\sum_{j=1}^n\sum_{k\neq j} h(X_j,X_k)\quad  \mbox
{and}\quad  n V_n=\frac{1}{n}\sum_{j,k=1}^n h(X_j,X_k),
\]
with $h\dvt \R^d\times\R^d\to\R$ symmetric and $\int_{\R^d}
h(x,y)P_X(\mathrm{d}x)=0, \forall y\in\R^d,$ can be derived by using a
spectral decomposition of the kernel, $h(x,y)=\sum_{k=1}^\infty\lambda
_k \Phi_k(x)\Phi_k(y)$, which holds true in the $L_2$-sense. Here,
$(\Phi_k)_k$ denote orthonormal eigenfunctions and $(\lambda_k)_k$ the
corresponding eigenvalues of the integral equation
%
\begin{equation}\label{eq.inteq}
\int_{\R^d}h(x,y)g(y)P_X(\mathrm{d}y)=\lambda  g(x).
\end{equation}
Approximate $n U_n$ by $n U_n^{(K)}=\sum_{k=1}^K \lambda_k \{
( n^{-1/2} \sum_{i=1}^n \Phi_k(X_i) )^2
 -  n^{-1} \sum_{i=1}^n \Phi_k^2(X_i) \}$. Then the sum under the
round brackets is asymptotically standard normal while the latter sum
converges in probability to 1. Finally, one obtains
%
\begin{equation}\label{eq.ustat}
n U_n
\stackrel{d}{\longrightarrow}\sum_{k=1}^\infty\lambda_k (Z_k^2-1),
\end{equation}
where $(Z_k)_{k}$ is a sequence of i.i.d.~standard normal random
variables; cf.~Serfling \cite{serf80}. If additionally $\E|h(X_1,X_1)|<\infty$,
the weak law of large numbers and Slutsky's theorem imply $V_n\stackrel
{d}{\longrightarrow}\sum_{k=1}^\infty\lambda_k (Z_k^2 -1)+ \E
h(X_1,X_1)$. (Here, $\stackrel{d}{\longrightarrow}$ denotes convergence
in distribution.)

So far, most previous attempts to derive the limit distributions of
degenerate $U$- and $V$-statistics of dependent random variables are
based on the adoption of this method of proof.
Eagleson \cite{E79} developed the asymptotic theory in the case of a strictly
stationary sequence of $\phi$-mixing, real-valued random variables
under the assumption of absolutely summable eigenvalues. This condition
is satisfied if the kernel function is of the form $h(x,y)=\int_\R
h_1(x,z)h_1(z,y)P_X(\mathrm{d}z)$ and $h_1$ is squared integrable w.r.t. $P_X$.
Using general heavy-tailed weight functions instead of $P_X$, the
eigenvalues are not necessarily absolutely summable; see, for example,
de Wet \cite{dew87}.
Carlstein \cite{Car88} analysed $U$-statistics of $\alpha$-mixing, real-valued
random variables in the case of finitely many eigenfunctions. He
derived a~limit distribution of the form (\ref{eq.ustat}), where
$(Z_k)_{k\in\N}$ is a sequence of centered normal random variables.
Denker \cite{Den82} considered stationary sequences $(X_n=f(Y_n,Y_{n+1},\ldots
))_n$ of functionals of $\beta$-mixing random variables $(Y_n)_n$. He
assumed $f$ and the cumulative distribution function of $X_1$ to be H\"
older continuous. Imposing some smoothness condition on~$h$, the limit\vadjust{\goodbreak}
distribution of $n U_n$ was derived under the additional
assumption~$\|\Phi_k\|_\infty<\infty$, $\forall  k\in\N$. The condition on $(\Phi
_k)_k$ is difficult or even impossible to check in a multitude of cases
since this requires to solve the associated integral equation~(\ref
{eq.inteq}). Similar difficulties occur if one wants to apply the
results of Dewan and Prakasa~Rao~\cite{DPR01} or Huang and Zhang \cite{HuZh06}. They studied $U$-statistics
of associated, real-valued random variables. Besides the absolute
summability of the eigenvalues, certain regularity conditions have to
be satisfied uniformly by the eigenfunctions in order to obtain the
asymptotic distribution of~$n U_n$.

A different approach was used by Babbel \cite{Ba89} to determine the limit
distribution of $U$-statistics of $\phi$- and $\beta$-mixing random
variables. She deduced the limit distribution via a Haar wavelet
decomposition of the kernel and empirical process theory without
imposing the critical conditions mentioned above. However, she presumed
that $\iint h(x,y) P_{X_k,X_{k+n}}(\mathrm{d}x, \mathrm{d}y)=0,  \forall  k\in\Z,
n\in\N$. This assumption does in general not hold true within our
applications in Section~\ref{S3}.
Moreover, this approach is not suitable when dealing with
$U$-statistics of $\tau$-dependent random variables since Lipschitz
continuity will be the crucial property of the (approximating) kernel
in order to exploit the underlying dependence structure.

\subsection{Main results}\label{SS22}

Let $(X_n)_{n\in\N}$ be a sequence of $\R^d$-valued random variables on
some probability space $(\Omega,{\mathcal A}, P)$ with common
distribution~$P_X$.
In this subsection, we derive the limit distributions of
\[
n U_n=\frac{1}{n}\sum_{j=1}^n\sum_{k\neq j} h(X_j,X_k) \quad \mbox
{and} \quad n V_n=\frac{1}{n}\sum_{j,k=1}^n h(X_j,X_k),
\]
where $h:\R^d\times\R^d\to\R$ is a symmetric function with $\int_{\R
^d} h(x,y)P_X(\mathrm{d}x)=0, \forall  y\in\R^d$.
In order to describe the dependence structure of $(X_n)_{n\in\N}$, we
recall the definition of the $\tau$-dependence coefficient for $\R
^d$-valued random variables of Dedecker and Prieur \cite{DP05}.

\begin{defn}\label{def1}
Let $(\Omega,{\mathcal A}, P)$ be a probability space, ${\mathcal M}$~a
sub-$\sigma$-algebra of ${\mathcal A}$ and $X$ an $\R^d$-valued random
variable. Assume that $\E\|X\|_{l_1} < \infty$, where $\|x\|_{l_1}=\sum
_{i=1}^d |x_i|,$ and define
\[
\tau({\mathcal M}, X)=\E\biggl( \sup_{f\in\Lambda_1(\R^d)} \biggl|\int
_{\R^d} f(x) P_{X|{\mathcal M}}(\mathrm{d}x)-\int_{\R^{d}} f(x)P_X(\mathrm{d}x)
\biggr|\biggr).
\]
Here, $P_{X|{\mathcal M}}$ denotes the conditional distribution of $X$
given ${\mathcal M} $ and $ \Lambda_1(\R^d)$ denotes the set of
1-Lipschitz functions from $\R^d$ to $\R$.
\end{defn}

We assume

\begin{enumerate}[(A1)]
\item[(A1)]
\begin{enumerate}[(ii)]
\item[(i)]
$(X_n)_{n\in\N}$ is a (strictly) stationary sequence of $\R^d$-valued
random variables on some probability space $(\Omega,{\mathcal A}, P)$
with common distribution $P_X$ and $\E\|X_1\|_{l_1} < \infty$.
\item[(ii)]
The sequence $(\tau_r)_{r\in\N}$, defined by
\begin{eqnarray*}
\tau_r &=& \sup\{\tau(\sigma(X_{s_1},\ldots,X_{s_u}),(X_{t_1}^\prime
,X_{t_2}^\prime,X_{t_3}^\prime)^\prime) |\\
&&\phantom{\sup\{}{} u\in\N, s_1\leq\cdots\leq s_u< s_u+r\leq t_1\leq t_2\leq t_3 \in\N\},
\end{eqnarray*}
satisfies $\sum_{r=1}^\infty r \tau_r^{\delta} <\infty$ for some
$\delta\in(0,1)$.
(Here, prime denotes the transposition.)
\end{enumerate}
\end{enumerate}

\begin{rem}\label{r.1}
If $\Omega$ is rich enough, due to Dedecker and Prieur \cite{DP04} the validity of (A1)
allows for the construction of a random vector $(\widetilde
X_{t_1}^\prime,\widetilde X_{t_2}^\prime,\widetilde X_{t_3}^\prime
)^\prime\stackrel{d}{=}(X_{t_1}^\prime,X_{t_2}^\prime
,X_{t_3}^\prime)^\prime$ that is independent of $X_{s_1},\ldots
,X_{s_u}$ and such that
%
\begin{equation}\label{eq.a1}
\sum_{i=1}^3 \E\|\widetilde X_{t_i}-X_{t_i}\|_{l_1}\leq\tau_r.
\end{equation}
\end{rem}

The notion of $\tau$-dependence is more general than mixing. If, for
example, $(X_n)_n$ is $\beta$-mixing, we obtain an upper bound for the
dependence coefficient $\tau_r\leq 6\int_0^{\beta(r)}Q_{|X_1|}(u)
\,\mathrm{d}u$, where $Q_{|X_1|}(u)=\inf\{t\in\R | P(\|X_1\|_{l_1}>t)\leq u\}
,  u\in[0,1],$ and $\beta(r)$ denotes the ordinary $\beta$-mixing
coefficient $\beta(r) :=\E\sup_{B\in\sigma(X_s, s\geq t+r), t\in\Z
}|P(B|\sigma(X_s, s\leq t))-P(B)|.$ This is a~consequence of Remark~2
of Dedecker and Prieur \cite{DP04}.
Moreover, inequality~(\ref{eq.a1}) immediately implies
%
\begin{equation}\label{eq.cov}
| \operatorname{cov}( h(X_{s_1},\ldots,X_{s_u}),k(X_{t_1},\ldots,X_{t_v})
) |\leq2\|h\|_\infty \operatorname{Lip}(k)\biggl\lceil\frac{v}{3}\biggr\rceil
\tau_r
\end{equation}
for $s_1\leq\cdots\leq s_u<s_u+r\leq t_1\leq\cdots\leq t_v\in\N$ and
for all functions $h\dvtx \R^u\to\R$ and $k\dvtx\allowbreak\R^v\to\R$ in ${\mathcal L}:=\{
f\dvtx \R^p\to\R\mbox{ for some } p\in\N |  \mbox{Lipschitz continuous
and bounded} \}$. Therefore, a~sequence of random variables that
satisfies (A1) is $((\tau_r)_r,\mathcal{L},\psi)$-weakly
dependent in the sense of Doukhan and Louhichi \cite{DL99} with $\psi(h,k,u,v)=2\|h\|
_\infty \operatorname{Lip}(k)\lceil\frac{v}{3}\rceil$. (Here and in the
sequel, $\operatorname{Lip}(g)$ denotes the Lipschitz constant of a generic function~$g$.)
A list of examples for $\tau$-dependent processes including causal
linear and functional autoregressive processes is provided by Dedecker and Prieur \cite{DP05}.

Besides the conditions on the dependence structure of $(X_n)_{n\in\N}$,
we make the following assumptions concerning the kernel:

\begin{enumerate}[(A2)]
\item[(A2)]
\begin{enumerate}[(ii)]
\item[(i)]
The kernel $h: \R^d\times\R^d \rightarrow\R$ is a symmetric,
measurable function and degenerate under $P_X$, that is, $\int_{\R^d}
h(x,y)P_X(\mathrm{d}x)=0, \forall y\in\R^d$.
\item[(ii)]
For a $\delta$ satisfying~(A1)(ii), the following moment constraints
hold true with some $\nu>(2-\delta)/(1-\delta)$ and an independent copy
$\widetilde X_1$ of $X_1$:
\[
\sup_{k\in\N}\E|h(X_1,X_{1+k})|^{\nu}<\infty \quad\mbox{and}\quad
\E|h(X_1,\widetilde X_1)|^{\nu}<\infty.
\]
\end{enumerate}
\item[(A3)] The kernel $h$ is Lipschitz continuous.
\end{enumerate}
Using an appropriate kernel truncation, it is possible to reduce the
problem of deriving the asymptotic distribution of $n U_n$ to
statistics with bounded kernel functions.

\begin{lem}\label{l.1}
Suppose that \textup{(A1)}, \textup{(A2)}, and \textup{(A3)} are fulfilled. Then there exists a
family of bounded functions $(h_c)_{c\in\R^+}$ satisfying \textup{(A2)} and \textup{(A3)}
uniformly such that
%
\begin{equation}\label{eq1}
\lim_{c\to\infty}\sup_{n\in\N} n^2 \E(U_n-U_{n,c})^2=0,
\end{equation}
where $U_{n,c}=n^{-2}\sum_{j=1}^n\sum_{k\neq j} h_c(X_j,X_k)$.
\end{lem}

After this simplification of the problem, we intend to develop a
decomposition of the kernel that allows for the application of a
central limit theorem (CLT) for weakly dependent random variables.
One could try to imitate the proof of the i.i.d.~case. According to the
discussion in the previous subsection, this leads to prerequisites that
can hardly be checked in numerous cases. Therefore, we do not use a
spectral decomposition of the kernel but a wavelet decomposition. It
turns out that Lipschitz continuity is the central property the kernel
function should satisfy in order to exploit~(\ref{eq.a1}). For this
reason, the choice of Haar wavelets, as they were employed by Babbel \cite
{Ba89}, is inappropriate in the present situation. Instead, the
application of Lipschitz continuous scale and wavelet functions is more
suitable.

In the sequel, let $\phi$ and $\psi$ denote scale and wavelet functions
associated with an one-dimensional multiresolution analysis. As
illustrated by Daubechies \cite{Dau02}, Section~8, these functions can be selected
in such a manner that they possess the following
properties:\looseness=1
\begin{enumerate}[(1)]
\item[(1)] $\phi$ and $\psi$ are Lipschitz continuous,
\item[(2)] $\phi$ and $\psi$ have compact support,
\item[(3)] $\int_{-\infty}^{\infty}\phi(x)\, \mathrm{d}x=1$ and $\int_{-\infty
}^{\infty}\psi(x) \,\mathrm{d}x=0.$
\end{enumerate}\looseness=0
It is well known that an orthonormal basis in $L_2(\R^d)$ can be
constructed from $\phi$ and~$\psi$. For this purpose, define $E:=\{0,1\}
^d\setminus\{0_d\}$, where $0_d$ denotes the $d$-dimensional null
vector. In addition, set
\[
\varphi^{(i)}:=
\cases{
\phi& \quad $\mbox{for } i=0,$\vspace*{2pt}\cr
\psi&\quad $\mbox{for } i=1$
}
\]
and define functions $\Psi^{(e)}_{j,k}\dvtx \R^d\to\R, j\in\Z
,k=(k_1,\ldots,k_d)^\prime\in\Z^d,$ by
\[
\Psi^{(e)}_{j,k}(x):=2^{jd/2}\prod_{i=1}^d \varphi
^{(e_i)}(2^{j}x_i-k_i)\qquad \forall e=(e_1,\ldots,e_d)^\prime\in E,
x=(x_1,\ldots,x_d)^\prime\in\R^d.
\]
The system
$(\Psi^{(e)}_{j,k})_{e \in E, j\in\Z, k\in\Z^d}$
is an orthonormal basis of $L_2(\R^d)$, see Wojtaszczyk \cite{Woj97}, Section~5. The
same holds true for
$
(\Phi_{0,k})_{k\in\Z^d }\cup(\Psi^{(e)}_{j,k})_{ j\geq
0,e \in E, k\in\Z^d},
$
where the functions $\Phi_{j,k}\dvtx \R^d\to\R$ are given by $\Phi
_{j,k}(x):=2^{jd/2}\prod_{i=1}^d \phi(2^{j}x_i-k_i), j\in\Z, k\in\Z^d$.

Now, an $L_2$-approximation of $n U_{n,c}$ by a statistic based on a
wavelet approximation of $h_c$ can be established. To this end, we
introduce $\widetilde h^{(K,L)}_c$ with
%
\begin{eqnarray}\label{eq.hckl}
\widetilde h^{(K,L)}_c(x,y)
&:=& \sum_{ k_1,k_2 \in\{-L,\ldots,L\}^d } \alpha^{(c)}_{k_1,k_2}\Phi
_{0,k_1}(x)\Phi_{0,k_2}(y)
\nonumber
\\[-8pt]
\\[-8pt]
\nonumber
&&{}+ \sum_{j=0}^{J(K)-1}\sum_{ k_1,k_2\in\{-L,\ldots,L\}^d}\sum_{e\in\bar
E} \beta_{j;k_1,k_2}^{(c,e)}\Psi_{j;k_1,k_2}^{(e)}(x,y),
\end{eqnarray}
where $\bar E:=(E\times E)\cup(E\times\{0_d\})\cup(\{0_d\}\times E)$,
\[
\Psi_{j;k_1,k_2}^{(e)}:=
\cases{
\Psi_{j,k_1}^{(e_1)} \Psi_{j,k_2}^{(e_2)} &\quad$\mbox{for } (e_1^\prime
,e_2^\prime)^\prime\in E\times E,$\vspace*{2pt}\cr
\Psi_{j,k_1}^{(e_1)} \Phi_{j,k_2} &\quad$\mbox{for } (e_1^\prime,e_2^\prime
)^\prime\in E\times\{0_d\},$\vspace*{2pt}\cr
\Phi_{j,k_1} \Psi_{j,k_2}^{(e_2)} &\quad$\mbox{for } (e_1^\prime,e_2^\prime
)^\prime\in\{0_d\}\times E,$
}
\]
$\alpha_{k_1,k_2}^{(c)}=\iint_{\R^d\times\R^d}h_c(x,y) \Phi
_{0,k_1}(x) \Phi_{0,k_2}(y) \,\mathrm{d}x \,\mathrm{d}y$ and $\beta
^{(c,e)}_{j;k_1,k_2}=\iint_{\R^d\times\R^d}h_c(x,y)\times  \Psi
_{j;k_1,k_2}^{(e)}(x,\allowbreak y) \,\mathrm{d}x \,\mathrm{d}y$.
We refer to the degenerate version of $\widetilde h_c^{(K,L)}$ as
$h_c^{(K,L)}$, given by
\begin{eqnarray*}
h_c^{(K,L)}(x,y)&:= & \widetilde h_c^{(K,L)}(x,y)-\int_{\R
^d}\widetilde h_c^{(K,L)}(x,y)P_X(\mathrm{d}x)-\int_{\R^d}\widetilde
h_c^{(K,L)}(x,y)P_X(\mathrm{d}y)\\
&&{}+\iint_{\R^d\times\R^d}\widetilde
h_c^{(K,L)}(x,y)P_X(\mathrm{d}x)P_X(\mathrm{d}y).
\end{eqnarray*}
The associated $U$-type statistic will be denoted by $U_{n,c}^{(K,L)}$.

\begin{lem}\label{l.5}
Assume that \textup{(A1)}, \textup{(A2}), and \textup{(A3)} are fulfilled. Then the sequence of
indices $(J(K))_{K\in\N}$ in (\ref{eq.hckl}) with $J(K)\longrightarrow
_{K\to\infty} \infty$ can be chosen such that
\[
\lim_{K\to\infty}\mathop{\lim\sup}_{L\to\infty}\sup_{n\in\N}n^2 \E
\bigl(U_{n,c}-U_{n,c}^{(K,L)}\bigr)^2= 0.
\]
\end{lem}

Employing the CLT of Neumann and Paparoditis \cite{NeuPa05} and the continuous mapping
theorem, we obtain the limit distribution of $n U_{n,c}^{(K,L)}$.
Finally, based on this result, the asymptotics of the $U$-type
statistic $n U_n$ can be derived. Moreover, a weak law of large
numbers (Lemma~\ref{l.lln} in Section~\ref{SS52}) allows for
deducing the limit distribution of $n V_n$ since $n V_n=n
U_n+n^{-1}\sum_{k=1}^n h(X_k,X_k)$.

Before stating the main result of this section, we introduce constants
$A_{k_1,k_2}:=\operatorname{cov}(\Phi_{0,k_1}(X_1), \Phi_{0,k_2}(X_1) )$ and
\[
B_{j;k_1,k_2}^{(c,e)}:=
\cases{
\operatorname{cov}\bigl(\Psi_{j,k_1}^{(e_1)}(X_1),\Psi_{j,k_2}^{(e_2)}(X_1)\bigr) &\quad$\mbox{for }
(e_1^\prime,e_2^\prime)^\prime\in E\times E,$\vspace*{2pt}\cr
\operatorname{cov}\bigl(\Psi_{j,k_1}^{(e_1)}(X_1),\Phi_{j,k_2}(X_1)\bigr) &\quad$\mbox{for }
(e_1^\prime,e_2^\prime)^\prime\in E\times\{0_d\},$\vspace*{2pt}\cr
\operatorname{cov}\bigl(\Phi_{j,k_1}(X_1),\Psi_{j,k_2}^{(e_2)}(X_1)\bigr) &\quad$\mbox{for }
(e_1^\prime,e_2^\prime)^\prime\in\{0_d\}\times E,$}
 \qquad j\in\Z,\  k_1,k_2\in\Z^d.
\]

\begin{thmm}\label{t.1}
Suppose that the assumptions \textup{(A1)}, \textup{(A2)}, and \textup{(A3)} are fulfilled. Then,
as $n\to\infty$,
\[
n U_n \stackrel{d}{\longrightarrow} Z
\]
with
\begin{eqnarray*}
Z&:=&\lim_{c\to\infty}\Biggl(
\sum_{ k_1,k_2 \in\Z^d}\alpha^{(c)}_{ k_1,k_2}[Z_{k_1} Z_{k_2}-A_{k_1,k_2}]
\\
&&\hphantom{\lim_{c\to\infty}(}{}+  \sum_{j = 0}^\infty
\sum_{k_1,k_2 \in\Z^d} \sum_{e=(e_1^\prime,e_2^\prime)^\prime\in\bar E}
\beta_{j; k_1,k_2}^{(c,e)} \bigl[Z_{j;k_1}^{(e_1)}
Z_{j;k_2}^{(e_2)}-B_{j;k_1,k_2}^{(c,e)}\bigr]\Biggr).
\end{eqnarray*}
Here, $(Z_k)_{k\in\Z^d}$ as well as $(Z_{j;k}^{(e)})_{j\geq0, k\in\Z
^d, e\in\{0,1\}^d}$ are centered and jointly normally distributed
random variables and the r.h.s.~converges in the $L_2$-sense.
If additionally $\E|h(X_1,X_1)|<\infty$, then
\[
n V_n \stackrel{d}{\longrightarrow} Z+\E h(X_1,X_1).
\]
\end{thmm}

As in the case of i.i.d.~random variables, the limit distributions of
$n U_n$ and $n V_n$ are, up to a constant, weighted sums of products
of centered normal random variables. In contrast to many other results
in the literature, the prerequisites of this theorem, namely moment
constraints and Lipschitz continuity of the kernel, can be checked
fairly easily in many cases. Nevertheless, the asymptotic distribution
has a complicated structure. Hence, quantiles can hardly be determined
on the basis of the previous result. However, we show in the following
section that the conditional distributions of the bootstrap
counterparts of $n U_n$ and $n V_n$, given $X_1,\ldots,X_n$, converge
to the same limits in probability.

Of course, the assumption of Lipschitz continuous kernels is rather
restrictive. Thus, we extend our theory to a more general class of
kernel functions. The costs for enlarging the class of feasible kernels
are additional moment constraints.

Besides (A1) and (A2), we assume
\begin{enumerate}[(A4)]
\item[(A4)]
\begin{enumerate}[(ii)]
\item[(i)] The kernel function satisfies
\[
|h(x,y)-h(\bar x,\bar y)|\leq f(x,\bar x,y,\bar y)[ \|x-\bar x\|
_{l_1}+\|y-\bar y\|_{l_1}]  \qquad \forall  x,\bar x,y,\bar y\in\R^d,
\]
where $f\dvtx\R^{4d}\to\R$ is continuous.
Moreover,
\[
\sup_{Y_1,\ldots,Y_5\sim P_X} \E\Bigl(\max_{a_1,a_2\in[-A,A]^d}
[f(Y_{1},Y_{2}+a_1,Y_{3},Y_{4}+a_2)]^\eta\|Y_{5}\|_{l_1}\Bigr)<\infty
\]
for $\eta:=1/(1-\delta)$ with $\delta$ satisfying (A2) and some $A>0$.
\item[(ii)]
$\sum_{r=1}^\infty r (\tau_r)^{\delta^2}<\infty.$
\end{enumerate}
\end{enumerate}
Even though the assumption (A4)(i) has a rather technical structure, it
is satisfied for example,~by polynomial kernel functions as long as the
sample variables have sufficiently many finite moments.
Analogous to Lemma~\ref{l.1} and Lemma~\ref{l.5}, the following
assertion holds.

\begin{lem}\label{l.4}
Suppose that \textup{(A1)}, \textup{(A2)}, and \textup{(A4)} are fulfilled. Then a family of
bounded kernels $(h_c)_c$ satisfying \textup{(A2)} and \textup{(A4)} uniformly and the
sequence of indices $(J(K))_{K\in\N}$ in (\ref{eq.hckl}) with
$J(K)\longrightarrow_{K\to\infty} \infty$ can be chosen such that
\[
\lim_{c\to\infty} \limsup_{K\to\infty} \limsup_{L\to\infty} \sup
_{n\in\N}  \E\bigl(U_{n}-U_{n,c}^{(K,L)}\bigr)^2=0.
\]
\end{lem}

This auxiliary result implies the analogue of Theorem~\ref{t.1} for
non-Lipschitz kernels.

\begin{thmm}\label{t.2}
Assume that \textup{(A1)}, \textup{(A2)}, and \textup{(A4)} are satisfied. Then, as $n\to\infty$,
\[
n U_n\stackrel{d}{\longrightarrow}Z,
\]
where $Z$ is defined as in Theorem~\ref{t.1}. If additionally $\E
|h(X_1,X_1)|<\infty$, then
\[
n V_n \stackrel{d}{\longrightarrow} Z+\E h(X_1,X_1).
\]
\end{thmm}

\section{Consistency of general bootstrap methods}\label{S3}

As we have seen in the previous section, the limit distributions of
degenerate $U$- and $V$-statistics have a rather complicated structure.
Therefore, in the majority of cases it is quite difficult to determine
quantiles, which are required in order to derive asymptotic critical
values of $U$- and $V$-type test statistics. The bootstrap offers a
suitable way of approximating these quantities.

Given $X_1,\ldots,X_n$, let $X^*$ and $Y^*$ denote vectors of bootstrap
random variables with values in $\R^{d_1}$ and $\R^{d_2}$.
In order to describe the dependence structure of the bootstrap sample,
we introduce, in analogy to Definition~\ref{def1},
\[
\tau^*(Y^*, X^*,x_n):=\E\biggl( \sup_{f\in\Lambda_1(\R^{d_1})}
\biggl|\int_{\R^{d_1}} f(x) P_{X^*|Y^*}(\mathrm{d}x)-\int_{\R^{d_1}}
f(x)P_{X^*}(\mathrm{d}x)\biggr|\big| \X_n=x_n\biggr)
\]
provided that $\E(\|X^*\|_{l_1}\vert\X_n=x_n)<\infty$ with $\X
_n:=(X_1^\prime,\ldots,X_n^\prime)^\prime$. We make the following assumptions:

\begin{enumerate}[$\mathrm{(A1^*)}$]
\item[$\mathrm{(A1^*)}$]
\begin{enumerate}[(ii)]
\item[(i)]The sequence of bootstrap variables is stationary with probability
tending to one. Additionally, $ (X^{*\prime}_{t_1},X^{*\prime
}_{t_2})^\prime\stackrel{d}{\longrightarrow} (X_{t_1}^\prime
,X_{t_2}^\prime)^\prime, \forall t_1,t_2 \in\N,$ holds true in probability.
\item[(ii)]Conditionally on $X_1,\ldots,X_n$, the random variables $(X_k^*)_{k\in\Z
}$ are $\tau$-weakly dependent, that is,~there exist a sequence of
coefficients $(\bar\tau_r)_{r\in\N}$ with $\sum_{r=1}^\infty r(\bar\tau
_r)^\delta<\infty$ for some $\delta\in(0,1)$, a constant $C_1<\infty$,
and a sequence~of sets $(\XX_n^{(1)})_{n\in\N}$ with $P(\X_n\in\XX
_n^{(1)})\longrightarrow_{n\to\infty} 1$ and the following property:
For any sequence $(x_n)_{n\in\N}$ with $x_n \in\XX_n^{(1)}, n\in\N$,
\mbox{$\sup_{k\in\N}\E(\|X_k^*\|_{l_1}\vert\X_n=x_n)\leq C_1$} and
\begin{eqnarray*}
\tau_r^*(x_n) &:=& \sup\{\tau^*((X^{*\prime}_{s_1},\ldots
,X^{*\prime}_{s_u})^\prime,(X_{t_1}^{*\prime},X_{t_2}^{*\prime
},X_{t_3}^{*\prime})^\prime, x_n) |\\[-2pt]
&&\hphantom{\sup\{}{} u\in\N,  s_1\leq\cdots\leq s_u< s_u+r\leq t_1\leq t_2\leq t_3 \in\N
\}
\end{eqnarray*}
can be bounded by $\bar\tau_r$ for all $r\in\N$.\vspace*{-1pt}
\end{enumerate}
\end{enumerate}

\begin{rem}
\begin{enumerate}[(ii)]
\item[(i)]
Neumann and Paparoditis \cite{NeuPa05} proved that in case of stationary Markov chains of
finite order, the key for convergence of the finite-dimensional
distributions is convergence of the conditional distributions, cf.
their Lemma~4.2. In particular, they showed that $\operatorname{AR}(p)$~bootstrap and
$\operatorname{ARCH}(p$)~bootstrap yield samples that satisfy (A$1^*$)(i).
\item[(ii)]
In Section~\ref{SS42}, we present another example that satisfies
(A$1^*$), namely a residual-based bootstrap procedure for a Lipschitz
contracting nonlinear $\operatorname{AR}(1)$~process, given by
$X_{t}=g(X_{t-1})+\varepsilon_t$. In particular, note that the
bootstrap process there cannot be proved to be mixing according to the
discreteness of the bootstrap innovations that are generated via
Efron's bootstrap from the empirical distribution of the recentered
residuals of the original process.\vspace*{-1pt}
\end{enumerate}
\end{rem}

\begin{lem}\label{l.7}
Suppose that \textup{(A1)} and \textup{(A$1^*$)} hold true. Further let \mbox{$h\dvtx \R^d\times\R
^d\to\R$} be a~bounded, symmetric, Lipschitz continuous function such
that
$\E h(X_1,y)=\E(h(X_1^*,y)|\allowbreak X_1,\ldots,X_n)= 0, \forall y\in\R^d$. Then,
\[
\frac{1}{n}\sum_{j=1}^n\sum_{k\neq j}h(X_j^*,X_k^*)\stackrel
{d}{\longrightarrow}Z\quad \mbox{and}\quad
\frac{1}{n}\sum_{j,k=1}^n h(X_j^*,X_k^*)\stackrel{d}{\longrightarrow
}Z+\E h(X_1,X_1)
\]
hold in probability as $n\to\infty$. Here, $Z$ is defined as in
Theorem~\textup{\ref{t.1}}.\vspace*{-1pt}
\end{lem}

In order to deduce bootstrap consistency, additionally, convergence in
a certain metric~$\rho$ is required, that is,
\[
\rho\Biggl(P\Biggl(\frac{1}{n}\sum_{j,k=1}^nh(X_j^*,X_k^*) \leq x
|X_1,\ldots,X_n\Biggr), P\Biggl(\frac{1}{n}\sum_{j,k=1}^nh(X_j,X_k) \leq x
\Biggr)\Biggr)\stackrel{P}{\longrightarrow}0.
\]
(Here, $\stackrel{P}{\longrightarrow}$ denotes convergence in
probability.) Convergence in the uniform metric follows from Lemma~\ref
{l.7} if the limit distribution has a continuous cumulative
distribution function. The next assertion gives a necessary and
sufficient condition for this.\vspace*{-1pt}

\begin{lem}\label{l.6}
The limit variable $Z$, derived in Theorem~\textup{\ref{t.1}}/Theorem~\textup{\ref
{t.2}} under \textup{(A1)}, \textup{(A2)}, and \textup{(A3)}/\textup{(A4)}, has a continuous cumulative
distribution function if $\operatorname{var}(Z)>0$.\vspace*{-1pt}
\end{lem}

Kernels of statistics emerging from goodness-of-fit tests for composite
hypotheses often depend on an unknown\vadjust{\goodbreak} parameter. We establish bootstrap
consistency for this setting, that is,~when parameters have to be
estimated. Moreover, the class of feasible kernels is enlarged.
For this purpose, we additionally assume

\begin{enumerate}[$\mathrm{(A2^*)}$]
\item[$\mathrm{(A2^*)}$]
\begin{enumerate}[(iii)]
\item[(i)] $\widehat\theta_n\stackrel{P}{\longrightarrow}\theta
\in\Theta\subseteq\R^p.$
\item[(ii)]
$\E(h(X_1^*,y,\widehat\theta_n)| \X_n )=0, \forall y\in\R^d$.
\item[(iii)]
For some $\delta$ satisfying (A$1^*$)(ii), $\nu>(2-\delta)/(1-\delta)$,
and a constant $C_2<\infty$, there exists a sequence of sets $(\XX
_n^{(2)})_{n\in\N}$ such that $P(\X_n\in\XX_n^{(2)})\longrightarrow
_{n\to\infty} 1$ and $\forall (x_n)_{n\in\N}$ with $x_n\in\XX
_n^{(2)}$ the following moment constraint holds true:
\[
\sup_{1 \leq k< n}\E\bigl(|h(X_1^*,X_{1+k}^*,\widehat\theta_n)|^{\nu
}+|h(X_1^*,\widetilde X_1^{*},\widehat\theta_n)|^{\nu}|\X_n=x_n
\bigr)\leq C_2,
\]
where (conditionally on $\X_n$) $\widetilde X_1^{*}$ denotes an
independent copy of $X_1^*$.
\end{enumerate}
\item[$\mathrm{(A3^*)}$]
\begin{enumerate}[(iii)]
\item[(i)] The kernel is continuous in its third argument in
some neighbourhood $U(\theta)\subseteq\Theta$ of $\theta$ and satisfies
\[
|h(x,y,\widehat\theta_n)-h(\bar x,\bar y,\widehat\theta_n)|\leq
f(x,\bar x,y,\bar y,\widehat\theta_n)[ \|x-\bar x\|_{l_1}+\|y-\bar
y\|_{l_1}]
\]
for all $x,\bar x,y,\bar y\in\R^d$, where $f\dvtx\R^{4d}\times\R^p\to\R$
is continuous on $\R^{4d}\times U(\theta)$.
Moreover, for $\eta:=1/(1-\delta)$ and some constants $A>0, C_3<\infty$
there exists a sequence of sets $(\XX_n^{(3)})_{n\in\N}$ such that
$P(\X_n\in\XX_n^{(3)})\longrightarrow_{n\to\infty} 1$ and $\forall
(x_n)_{n\in\N}$ with $x_n\in\XX_n^{(3)}$ the following moment
constraint holds true:\looseness=1
\[
\E\Bigl(\max_{a_1,a_2\in[-A,A]^d}
[f(Y_{1}^*,Y_{2}^*+a_1,Y_{3}^*,Y_{4}^*+a_2,\widehat\theta_n)]^\eta\|
Y_{5}^*\|_{l_1}\big|\X_n=x_n\Bigr)\leq C_3
\]\looseness=0
for all $Y_1^*,\ldots, Y_5^*$ with $Y_k^*\stackrel{d}{=}X_1^*, k\!\in\!\{
1,\ldots,5\}$ (conditionally on $X_1,\ldots,X_n$).
\item[(ii)]
$\sum_{r=1}^\infty r(\bar\tau_r)^{\delta^2} <\infty$.
\end{enumerate}
\end{enumerate}
Under these assumptions a result concerning the asymptotic
distributions of $n U_n^*=n^{-1}\times \sum_{j=1}^n\sum_{k\neq j}
h(X_j^*,X_k^*, \widehat\theta_n)$ and $n V_n^*=n^{-1}\sum_{j,k=1}^n
h(X_j^*,X_k^*, \widehat\theta_n)$ can be derived. To this end, we
denote the $U$- and $V$-statistics with kernel $h(\cdot,\cdot,\theta)$
and arguments $X_1,\ldots, X_n$ by $U_n$ and $V_n$, respectively.

\begin{thmm}\label{t.3}
Suppose that the conditions \textup{(A$1$)}, \textup{(A$2$)}, and \textup{(A$4$)} as well as
\textup{(A$1^*$)}, \textup{(A$2^*$)}, and \textup{(A$3^*$)} are fulfilled.
\begin{enumerate}[(ii)]
\item[(i)] As $n\to\infty$,
\[
n U_n^* \stackrel{d}{\longrightarrow} Z,\qquad   \mbox{in probability,}
\]
where $Z$ is defined as in Theorem~\ref{t.1}. If furthermore
$\operatorname{var}(Z)>0$, then
\[
\sup_{-\infty< x < \infty}|P(n U_n^* \leq x |X_1,\ldots,X_n)- P(n
U_n \leq x)| \stackrel{P}{\longrightarrow}0.
\]
\item[(ii)]
If additionally $\E|h(X_1,X_1,\theta)|<\infty$ and
$
\E(|h(X_1^*,X_1^*,\widehat\theta_n)| \vert\X_n )\stackrel
{P}{\longrightarrow}\E|h(X_1,X_1,\theta)|,
$
then as $n\to\infty$,
\[
n V_n^* \stackrel{d}{\longrightarrow} Z+ \E h(X_1,X_1,\theta),\qquad
\mbox{in probability}.
\]
Moreover, in case of $\operatorname{var}(Z)>0$,
\[
\sup_{-\infty< x < \infty}|P(n V_n^* \leq x |X_1,\ldots,X_n)- P(n
V_n \leq x)| \stackrel{P}{\longrightarrow}0.
\]
\end{enumerate}
\end{thmm}

\begin{rem}
Theorem~\ref{t.3} implies that bootstrap-based tests of $U$- or
$V$-type have asymptotically a prescribed size $\alpha$, that is, $P(n
U_n>t_{u,\alpha}^*)\longrightarrow_{n\to\infty} \alpha$ and $P(n
V_n>t_{v,\alpha}^*)\longrightarrow_{n\to\infty} \alpha$, where
$t^*_{u,\alpha}$ and $t^*_{v,\alpha}$ denote the $(1-\alpha)$-quantiles
of $n U_n^*$ and $n V_n^*$, respectively, given $X_1,\ldots,X_n$.
\end{rem}

\section{$L_2$-tests for weakly dependent observations}\label{S4}

This section is dedicated to two applications in the field of
hypothesis testing. For sake of simplicity, we restrict ourselves to
real-valued random variables and consider simple null hypotheses only.
The test for symmetry as well as the model-specification test can be
extended to problems with composite hypotheses, cf. Leucht \cite{Leu10a,Leu10b}.

\subsection{A test for symmetry}\label{SS41}

Answering the question whether a distribution is symmetric or not is
interesting for several reasons. Often robust estimators of and
robust tests for location parameters assume the observations to arise
from a symmetric
distribution, see, for example, Staudte and Sheather \cite{StSh90}. Consequently, it is
important to check
this assumption before applying those methods. Moreover, symmetry
plays a central role in analyzing and modeling real-life phenomena. For
instance, it is often presumed that an observed process can be
described by an $\operatorname{AR}(p)$~process with Gaussian innovations which in turn
implies a~Gaussian marginal distribution. Rejecting the hypothesis of
symmetry contradicts this type of marginal distribution. Furthermore,
this result of the test excludes any kind of symmetric innovations in
that context.

Suppose that we observe $X_1,\ldots, X_n$ from a sequence of real-valued
random variables with common distribution~$P_X$ and satisfying~(A1).
For some $\mu\in\R$, we are given the problem
\[
{\mathcal H}_0\dvt P_{X-\mu}=P_{\mu-X}
\quad\mbox{vs.}\quad {\mathcal H}_1\dvt  P_{X-\mu}\neq P_{\mu-X}.
\]
Similar to Feuerverger and Mureika \cite{FM77}, who studied the problem for i.i.d.~random
variables, we propose the following test statistic:
\[
S_n=n \int_\R\bigl[\Im\bigl(c_n(t)\mathrm{e}^{-\mathrm{i}\mu t}\bigr)\bigr]^2  w(t) \,\mathrm{d}t
=\frac{1}{n}\sum_{j,k=1}^n \int_\R\sin\bigl(t(X_j-\mu)\bigr) \sin\bigl(t(X_k-\mu)\bigr)
w(t)\, \mathrm{d}t
\]
which makes use of the fact that symmetry of a distribution is
equivalent to a vanishing imaginary part of the associated
characteristic function. Here, $\Im(z)$ denotes the imaginary part of
$z\in\C$, $c_n$ denotes the empirical characteristic function and $w$
is some positive measurable weight function with $\int_\R(1+|t|) w(t)
\,\mathrm{d}t <\infty.$ Obviously, $S_n$ is a $V$-type statistic whose kernel
satisfies (A2) and~(A3). Thus, its limit distribution can be determined
by Theorem~\ref{t.1}. Assuming that the observations come from a
stationary $\operatorname{AR}(p)$ or $\operatorname{ARCH}(p)$~process, the validity of (A$1^*$) is
assured when the $\operatorname{AR}(p)$ or $\operatorname{ARCH}(p)$~bootstrap methods given by Neumann and Paparoditis \cite
{NeuPa05} are used in order to generate the bootstrap counterpart of
the sample. Hence, in these cases the prerequisites of Lemma~\ref{l.7}
are satisfied excluding degeneracy. Inspired by Dehling and Mikosch~\cite{DM94}, who
discussed this problem for Efron's Bootstrap in the i.i.d.~case, we
propose a bootstrap statistic with the kernel
\[
h_n^*(x,y)=h(x,y)-\int_\R h(x,y)P_n^*(\mathrm{d}x)-\int_\R h(x,y)P_n^*(\mathrm{d}y)
+\int_{\R^2} h(x,y)P_n^*(\mathrm{d}x)P_n^*(\mathrm{d}y).
\]
Here, $h$ denotes the kernel function of $S_n$ and $P_n^*$ the
distribution of $X_1^*$ conditionally on $X_1,\ldots,X_n$.
Similar to the proof of Theorem~\ref{t.3}, the desired convergence
property of~$S_n^*$ can be verified.

\subsection{A model-specification test}\label{SS42}

Let $X_0,\ldots,X_n$ be observations resulting from a stationary
real-valued nonlinear autoregressive process with centered
i.i.d.~innovations~$(\varepsilon_k)_{k\in\Z}, $ that
is,~$X_k=g(X_{k-1})+\varepsilon_k.$ Suppose that $\E|\varepsilon
_0|^{4+\delta}\,{<}\,\infty$ for some $\delta\,{>}\,0$ and that
$g\,{\in}\, G\,{:=}\,\{f\dvtx \R\,{\to}\,\R |  f \mbox{ Lipschitz continuous}$
$\mbox{with } \operatorname{Lip}(f)<1\}$. Thus, the process $(X_k)_{k\in\Z}$ is $\tau
$-dependent with exponential rate, see Dedecker and Prieur \cite{DP05}, Example~4.2.
We will present a test for the problem
\[
{\mathcal H}_0\dvt  P\bigl(\E(X_1|X_{0})=g_0(X_{0})\bigr)=1
\quad\mbox{vs.}\quad {\mathcal H}_1\dvt  P\bigl(\E
(X_1|X_{0})=g_0(X_{0})\bigr)<1
\]
with $g_0\in G$.
For sake of simplicity, we stick to these small classes of functions
$G$ and of processes $(X_k)_{k\in\Z}.$ An extension to a more
comprehensive variety of model-specification tests is investigated in a
forthcoming paper, cf.~Leucht \cite{Leu10b}.

Similar to Fan and Li \cite{FL99}, we propose the following test statistic:
\begin{eqnarray*}
T_n&=&\frac{1}{n \sqrt h}\sum_{j=1}^n\sum_{k\neq
j}\bigl(X_j-g_0(X_{j-1})\bigr)\bigl(X_k-g_0(X_{k-1})\bigr)K\biggl(\frac
{X_{j-1}-X_{k-1}}{h}\biggr)\\
&=:& \frac{1}{n}\sum_{j=1}^n\sum_{k\neq j} H(Z_j,Z_k),
\end{eqnarray*}
that is,~a kernel estimator (multiplied with $n\sqrt h$) of $\E
([X_1-g(X_0)]\E(X_1-g(X_0)|\break X_0)p(X_0))$ that is equal to zero under
${\mathcal H_0}$. Here, $Z_k:=(X_k,X_{k-1})^\prime, k\in\Z,$ and $p$
denotes the density of the distribution of $X_0$.\vadjust{\goodbreak}
Fan and Li \cite{FL99}, who considered $\beta$-mixing processes, used a similar
test statistic with a vanishing bandwidth. In contrast, we consider the
case of a fixed bandwidth. These tests are more powerful against Pitman
alternatives $g_{1,n}(x)=g_0(x)+n^{-\beta}w(x)+\mathrm{o}(n^{-\beta}), \beta
>0, w\in G$. For a detailed discussion of this topic, see Fan and Li \cite{FL00}.

Obviously, $T_n$ is degenerate under ${\mathcal H_0}$. If we assume $K$
to be a bounded, even, and Lip\-schitz continuous function, then there
exists a function $f\dvtx\R^8\to\R$ with $ |H(z_1,z_2)-H(\bar z_1,\bar
z_2)|\leq f(z_1,\bar z_1,z_2,\bar z_2)(\|z_1-\bar z_1\|_{l_1}+\|
z_2-\bar z_2\|_{l_1})$ and such that (A4) is valid. Moreover, under
these conditions $H$ satisfies (A2). Hence, the assertion of
Theorem~\ref{t.2} holds true.
In order to determine critical values of the test, we propose the
bootstrap procedure given by Franke and Wendel~\cite{FW92} (without estimating the
regression function). The bootstrap innovations $(\varepsilon_t^*)_{t}$
are drawn with replacement from the set $\{\tilde\varepsilon
_t=\varepsilon_t-n^{-1}\sum_{k=1}^n\varepsilon_k\}_{t=1}^n$, where
$\varepsilon_t=X_t-g_0(X_{t-1}),  t=1,\ldots,n$. After choosing a
starting value $X_0^*$ independently of $(\varepsilon_t^*)_{t\geq1}$,
the bootstrap sample $X_t^*=g(X_{t-1}^*)+\varepsilon_t^*$ as well as
the bootstrap counterpart $T_n^*= n^{-1}\sum_{j=1}^n\sum_{k\neq j}
H(Z_j^*,Z_k^*)$ of the test statistic with
$Z_k^*=(X_k^*,X_{k-1}^*)^\prime,  k=1,\ldots,n,$ can be computed. In
contrast to the previous subsection, the proposed bootstrap method
leads to a degenerate kernel function.
Obviously, the bootstrap sample is $\tau$-dependent in the sense of
(A$1^*$) and satisfies \mbox{$\E(|X_k^*| | Z_1,\ldots,Z_n)<C$} for some
$C<\infty$ with probability tending to one.
Theorem~1 of Diaconis and Freedman~\cite{DF99} yields the existence of a stationary solution
to $X_t^*=g(X_{t-1}^*)+\varepsilon_t^*$ and that the distribution of
any ``reasonably'' started process converges to the stationary one with
exponential rate. In order to apply our theory, $X_0^*$ is assumed to
be drawn from the stationary bootstrap distribution, conditionally on
$X_1,\ldots,X_n$.
We employ Lemma~4.2 of Neumann and Paparoditis \cite{NeuPa05} to verify convergence of the
finite dimensional distributions. The application of this result
requires the convergence of the conditional distributions, that is, ~$
\sup_{x\in K} d(P^{X_t^*|X_{t-1}^*=x},P^{X_t|X_{t-1}=x})\stackrel
{P}{\longrightarrow}0$ for every compact $K\subset\R$ and $d(P,Q)=\inf
_{X\sim P,Y\sim Q}\E(|X-Y|\wedge1).$ In the present context, this can
be confirmed similarly to the proof of Lemma~4.1 by Neumann and Paparoditis \cite{NeuPa05} if
the innovations of the original process have a bounded density. Summing
up, all prerequisites of Theorem~\ref{t.3} are satisfied. Hence,
critical values of the above test can be determined using the proposed
model-based bootstrap procedure.\vspace*{-2pt}

\section{Proofs}\label{S5}\vspace*{-2pt}

\subsection{Proofs of the main theorems}\label{SS51}\vspace*{-2pt}

Throughout this section, $C$ denotes a positive finite generic constant.

\begin{pf*}{Proof of Theorem~\protect\ref{t.1}}
First, we derive the limit distribution of $n U_{n,c}^{(K,L)}$,
defined before Lemma~\ref{l.5}. Afterwards, the asymptotic
distributions of $n U_n$ and $n V_n$ are deduced by means of
Lemma~\ref{l.1}, Lemma~\ref{l.5}, and a weak law of large numbers.

The following modified representation of $\widetilde h_c^{(K,L)}$ will
be useful in the sequel:
\[
\widetilde h^{(K,L)}_c(x,y) = \sum_{k,l=1}^{M(K,L)} \gamma_{k,l}^{(c)}
\tilde q_k(x) \tilde q_l(y),\vadjust{\goodbreak}
\]
where $(\tilde q_l)_{l=1}^{M(K,L)}$ is an ordering of $\bigcup_{k\in\{
-L,\ldots,L\}^d}\{\{\Phi_{j,k}\} \cup\{\Psi^{(e)}_{j,k}\}_{e\in
E,j\in\{0,\ldots,J(K)-1\}}\}$ and $\gamma_{k,l}^{(c)}=\gamma
_{l,k}^{(c)}, k,l\in\{1,\ldots,M(K,L)\}$, are the associated
coefficients. Moreover, the introduction of
$ q_k(X_i) := \tilde q_k(X_i) - \E\tilde q_k(X_i),  k\in\{1,\ldots
,M(K,L)\}, i\in\{1,\ldots,n\}$,
allows for the compact notation of $n U_{n,c}^{(K,L)}$,
\[
n U_{n,c}^{(K,L)}= \sum_{k,l=1}^{M(K,L)} \gamma_{k,l}^{(c)}\Biggl(
\Biggl[\frac{1}{\sqrt n}\sum_{i=1}^n q_k(X_i)\Biggr]
\Biggl[\frac{1}{\sqrt n}\sum_{j=1}^n q_l(X_j)\Biggr]-\frac{1}{n}\sum
_{i=1}^n q_k(X_i)q_l(X_i)\Biggr).
\]
The latter summand in the round brackets converges to $-\E
q_k(X_1)q_l(X_1)$ in probability by virtue of Lemma~\ref{l.lln}.
In order to derive the limit distributions of the first summands,
we consider
$ n^{-1/2}\sum_{i=1}^n (q_1(X_i),\ldots,q_{M(K,L)}(X_i))^\prime$.
Due to the Cram\'er--Wold device, it suffices to investigate $\sum
_{k=1}^{M(K,L)} t_k n^{-1/2}\sum_{i=1}^n q_k(X_i),$ $\forall
(t_1,\ldots,t_{M(K,L)})^\prime\in\R^{M(K,L)}$.
Asymptotic normality can be established by applying the CLT of Neumann and Paparoditis \cite
{NeuPa05} to $Q_i:= \sum_{k=1}^{M(K,L)} t_k q_k(X_i), i=1,\ldots,n$. To
this end, the prerequisites of this tool have to be checked. Obviously,
we are given a strictly stationary sequence of centered bounded random
variables. This implies in conjunction with the dominated convergence
theorem that the Lindeberg condition is fulfilled.
In order to show
\[
\frac{1}{n} \operatorname{var} (Q_1+ \cdots+ Q_n) \ninfty\sigma^2 := \operatorname{var}(Q_1)+ 2\sum
_{k=2}^\infty \operatorname{cov}(Q_1,Q_k),
\]
the validity of (A1) can be employed which moreover assures the
existence of the limit~$\sigma^2$. Then,
\begin{eqnarray*}
\biggl| \frac{1}{n} \operatorname{var}(Q_1+ \cdots + Q_n) - \sigma^2 \biggr|
&=& \Biggl| \frac{2}{n} \sum_{r=2}^{n} (n-[r-1]) \operatorname{cov}(Q_1,Q_r)- 2\sum
_{k=2}^\infty \operatorname{cov}(Q_1,Q_k) \Biggr| \\
&\leq& 2 \sum_{r=2}^\infty\min\biggl\{\frac{r-1}{n},1\biggr\} |
\operatorname{cov}(Q_1,Q_r)|\\
&\leq& 4\|Q_1\|_\infty \operatorname{Lip}(Q_1)\sum_{r=2}^\infty\min\biggl\{\frac
{r-1}{n},1\biggr\}\tau_{r-1},
\end{eqnarray*}
where the latter inequality follows from~(\ref{eq.cov}). The
summability condition of the dependence coefficients in connection with
Lebesgue's dominated convergence theorem yields the desired result.
Since $Q_{t_1}Q_{t_2}$ forms a Lip\-schitz continuous function,
inequality~(6.4) of Neumann and Paparoditis \cite{NeuPa05} holds true with $\theta
_r=\operatorname{Lip}(Q_{t_1}Q_{t_2}) \tau_r$. It is easy to convince oneself that
their condition~(6.3) is not needed if the involved random variables
are uniformly bounded. Finally, we obtain
\[
n^{-1/2}(Q_1+ \cdots + Q_n)\stackrel{d}{\longrightarrow}N(0,\sigma^2)
\]
and hence,
\begin{eqnarray*}
n  U_{n,c}^{(K,L)} &\stackrel{d}{\longrightarrow}& Z^{(K,L)}_c \\[-3pt]
&\hspace*{4pt}:=&
\sum_{k_1,k_2\in\{-L,\ldots,L\}^d}\alpha^{(c)}_{k_1,k_2}[Z_{k_1}
Z_{k_2}-A_{k_1,k_2}]\\[-3pt]
&&{}+  \sum_{j=0}^{J(K)-1}
\sum_{k_1,k_2\in\{-L,\ldots,L\}^d} \sum_{e=(e_1^\prime,e_2^\prime)^\prime
\in\bar E}
\beta_{j; k_1,k_2}^{(c,e)}\bigl[ Z_{j;k_1}^{(e_1)}
Z_{j;k_2}^{(e_2)}-B_{j;k_1,k_2}^{(e)}\bigr].\vspace*{-2pt}
\end{eqnarray*}
Here, $(Z_k)_{k\in\{-L,\ldots,L\}^d}$ and $(Z_{j;k}^{(e)})_{j\in\{
0,\ldots,J(K)-1\},e\in\{0,1\}^d,k\in\{-L,\ldots,L\}^d}$,\vspace*{1pt} respectively, are
centered and jointly normally distributed random variables.

By Lemma~\ref{l.1} and Lemma~\ref{l.5}, we have\vspace*{-1pt}
\[
\lim_{c\to\infty}\limsup_{K\to\infty}\limsup_{L\to\infty}\sup_{n\in\N
} n^2 \E\bigl(U_{n,c}^{(K,L)}-U_n\bigr)^2=0.\vspace*{-2pt}
\]
Since $n  U_{n,c}^{(K,L)} \stackrel{d}{\longrightarrow} Z^{(K,L)}_c $,
it remains to show
%
\begin{equation}\label{eq.Z}
\lim_{c\to\infty}\limsup_{K\to\infty}\limsup_{L\to\infty} \E
\bigl(Z_{c}^{(K,L)}-Z\bigr)^2=0\vspace*{-2pt}
\end{equation}
in order to prove that $n U_n\stackrel{d}{\longrightarrow} Z$ due to
Billingsley \cite{Bil68}, Theorem~4.2. To this end, we first show that
$(Z^{(K,L)}_c)_L$ is a Cauchy sequence in $L_2$.
Note that $n  (U_{n,c}^{(K,L_1)}-U_{n,c}^{(K,L_2)})\stackrel
{d}{\longrightarrow} Z^{(K,L_1)}_c-Z^{(K,L_2)}_c$.
According to Theorem~5.3 of Billingsley \cite{Bil68}, we obtain
$\E(Z^{(K,L_1)}_c-Z^{(K,L_2)}_c)^2 \leq\lim\inf_{n\to\infty} n^2 \E
(U_{n,c}^{(K,L_1)}-U_{n,c}^{(K,L_2)})^2.$ The r.h.s.~converges to zero
as $L_1,L_2\to\infty$ by virtue of~(\ref{eq.l5.2}) in the proof of
Lemma~\ref{l.5}.
Denoting the corresponding limit by~$Z_c^{(K)}$ similar arguments yield
\begin{eqnarray*}
\E\bigl(Z^{(K_1)}_c-Z_c^{(K_2)}\bigr)^2
&\leq& 4\limsup_{L\to\infty} \E\bigl(Z_c^{(K_1,L)} -Z_c^{(K_2,L)}
\bigr)^2\\[-2pt]
&\leq& 4\limsup_{L\to\infty} \liminf_{n\to\infty} n^2 \E
\bigl(U_{n,c}^{(K_1,L)} -U_{n,c}^{(K_2,L)}\bigr)^2\\[-2pt]
&\leq& 16\liminf_{n\to\infty}n^2 \E
\bigl(U_{n,c}^{(K_1)}-U_{n,c}^{(K_2)}\bigr)^2 \mathop{\longrightarrow
} _{K_1,K_2\to\infty}  0\vspace*{-2pt}
\end{eqnarray*}
according to (\ref{eq.l5.1}) of the proof of Lemma~\ref{l.5}.
In view of Lemma~\ref{l.1}, we obtain (\ref{eq.Z}) by applying the
above method once again.
This in turn leads to the desired limit distribution of $n U_n$.

Based on the result concerning $U$-type statistics, the limit
distribution of $n V_n$ can be established. Since
$V_n=U_n + n^{-2}\sum_{k=1}^{n} h(X_k,X_k)$, it remains to verify
that $n^{-1}\sum_{k=1}^{n} h(X_k,X_k)\stackrel{P}{\longrightarrow}\E
h(X_1,X_1).$ This in turn is a consequence of Lemma~\ref{l.lln}.\vspace*{-2pt}
\end{pf*}

\begin{pf*}{Proof of Theorem~\protect\ref{t.2}}
On the basis of Lemma~\ref{l.4} similar arguments as in the proof of
Theorem~\ref{t.1} yield $n U_n\stackrel{d}{\longrightarrow}Z$.
Moreover, Lemma~\ref{l.lln} implies $n^{-1}\sum_{k=1}^n
h(X_k,X_k)\stackrel{P}{\longrightarrow} \E h(X_1,X_1)$. Thus, $n
V_n\stackrel{d}{\longrightarrow}Z+\E h(X_1,X_1)$.\vadjust{\goodbreak}
\end{pf*}

\begin{pf*}{Proof of Theorem~\textup{\protect\ref{t.3}}}
Due to Lemma~\ref{l.6}, it suffices to verify distributional convergence.
To this end, we introduce
\[
\XX_n^{\theta} \subseteq\XX_n^{(1)}\cap\XX_n^{(2)}\cap\XX_n^{(3)}\cap
\{\X_n| \|\widehat\theta_n-\theta\|_{l_1}<\delta_n\}
\]
such that
%
\begin{eqnarray}
{\mathcal L}\bigl((X_{t_1}^{*\prime},\ldots,X_{t_k}^{*\prime})^\prime|\X_n=x_n\bigr)&=&
{\mathcal L}\bigl((X_{t_1+l}^{*\prime},\ldots,X_{t_k+l}^{*\prime})^\prime|\X
_n=x_n\bigr), \label{eqn.1}\\
{\mathcal L}\bigl((X_{t_1}^{*\prime},X_{t_2}^{*\prime})^\prime|\X
_n=x_n\bigr)&\Longrightarrow&
{\mathcal L}((X_{t_1}^{\prime},X_{t_2}^{\prime})^\prime)\label{eqn.2}
\end{eqnarray}
uniformly for any sequence $(x_n)_{n\in\N}$ with $x_n\in\XX_n^\theta$
and $t_1,\ldots,t_k,k,l\in\N$.
Moreover, the null sequence $(\delta_n)_{n\in\N}$ can be chosen such
that on $\XX_n^\theta,$ $\widehat\theta_n\in U(\theta)$ and $P(\X_n
\in\XX_n^{\theta})\longrightarrow_{n\to\infty} 1$ hold. Hence, to
prove $n U_n^*\stackrel{d}{\longrightarrow} Z$, in probability, it
suffices to verify that $n U_n^*$ converges to $Z$ in distribution
conditionally on $\X_n=x_n$ for any sequence $(x_n)_n$ with $x_n\in\XX
_n^\theta.$ Now, we take an arbitrary sequence $(x_n)_n$ with $x_n\in\XX
_n^\theta, n\in\N.$

In order to show that it suffices to investigate statistics with
bounded kernels, we consider the degenerate version $h^{*}_c$ of
\[
\widetilde h^{*}_c(x,y,\widehat\theta_n) :=
\cases{
h(x,y,\widehat\theta_n) &\quad$\mbox{for } |h(x,y,\widehat\theta_n)|\leq
c_h(\widehat\theta_n),$\vspace*{2pt}\cr
-c_h(\widehat\theta_n) &$\quad\mbox{for } h(x,y,\widehat\theta_n)
<-c_h(\widehat\theta_n),$\vspace*{2pt}\cr
c_h(\widehat\theta_n) &\quad$\mbox{for } h(x,y,\widehat\theta_n)
>c_h(\widehat\theta_n)$}
\]
with $c_h(\widehat\theta_n):=\max_{x,y\in[-c,c]^d} |h(x,y,\widehat\theta
_n)|\leq\max_{x,y\in[-c,c]^d,\|\bar\theta\|_{l_1}\leq\delta_1}
|h(x,y,\bar\theta)|<\infty$. The associated $U$-statistics are denoted
by $U^*_{n,c}$. Now, imitating the proof of Lemma~\ref{l.1} results in
\[
\mathop{\lim\sup}_{n\to\infty} n^2 \E[(U_n^*-U_{n,c}^*)^2|\X_n=x_n] \mathop
{\longrightarrow} _{c\to\infty}0.
\]
Within the calculations, the relation $\limsup_{n\to\infty}P
(X_1^*\notin (-c,c)^d\vert\X_n=x_n)\leq P(X_1\notin
(-c,c)^d)\longrightarrow_{c\to\infty}0$ has to be invoked which
follows from Portmanteau's theorem in conjunction with~(\ref{eqn.2}).
Next, we approximate the bounded kernel by the degenerate version of
\[
\widetilde h^{*(K,L)}_c
:= \sum_{ k_1,k_2 \in\{-L,\ldots,L\}^d } \widehat\alpha
^{(c)}_{k_1,k_2}\Phi_{0,k_1}\Phi_{0,k_2}+ \sum_{j=0}^{J(K)-1}\sum_{
k_1,k_2 \in\{-L,\ldots,L\}^d}\sum_{e\in\bar E} \widehat\beta
_{j;k_1,k_2}^{(c,e)}\Psi_{j;k_1,k_2}^{(e)},
\]
where $\widehat\alpha_{k_1,k_2}^{(c)}=\iint_{\R^d\times\R
^d}h^{*}_c(x,y, \widehat\theta_n) \Phi_{0,k_1}(x) \Phi_{0,k_2}(y)\, \mathrm{d}x\, \mathrm{d}y$
and $\widehat\beta^{(c,e)}_{j;k_1,k_2}=\iint_{\R^d\times\R
^d}h^{*}_c(x,y,\allowbreak\widehat\theta_n) \Psi_{j;k_1,k_2}^{(e)}(x,y)\, \mathrm{d}x\, \mathrm{d}y$.
Denoting the associated $U$-statistic by $\widehat U_{n,c}^{*(K,L)}$
leads to
\[
\lim_{K\to\infty}\limsup_{L\to\infty}\limsup_{n\to\infty}
n^2 \E\bigl[\bigl( U_{n,c}^*-\widehat U_{n,c}^{*(K,L)}\bigr)^2|\X_n=x_n\bigr]=0
\]
which can be proved by following the lines of the proof of Lemma~\ref{l.4}.
Here,~$J(K)$ is chosen as follows: We first select some $b=b(K)<\infty
$ such that $P(X_1\notin(-b,b)^d)\leq1/K$. Afterwards, we choose
$J(K)$ such that $\max_{x,y\in[-b,b]^d}|h_c(x,y,\theta)-\widetilde
h^{(K)}_c(x,y,\theta)|\leq{1/K}$ and $S_\phi/2^{J(K)}<A$, where $S_\phi
$ denotes the length of the support of the scale function $\phi$. The
index $J(K)$ can be determined independently of $n$ on $(\XX_n^\theta
)_n$ since $\max_{x,y\in[-b,b]^d}|h^{*}_c(x,y,\widehat\theta_n)-
h_c(x,y,\theta)|\longrightarrow0$ and $\max_{x,y\in[-b,b]^d}|\widetilde
h^{(K)}_c(x,y,\theta)-\widetilde h^{*(K)}_c(x,y,\allowbreak\widehat\theta
_n)|\longrightarrow0$, as $n\to\infty$, due to the continuity
assumptions on $f$. Here, $\widetilde h^{*(K)}_c$ is defined by the
substitution of $\sum_{ k_1,k_2 \in\{-L,\ldots,L\}^d }$ through $\sum_{
k_1,k_2 \in\Z^d }$ in the definition of $
\widetilde h^{*(K,L)}_c $.
Also note that
\begin{eqnarray*}
\widehat\alpha_{k_1,k_2}^{(c)}&\displaystyle\ninfty&\alpha_{k_1,k_2}^{(c)}:=\iint_{\R
^d\times\R^d}h_c(x,y,\theta) \Phi_{0,k_1}(x) \Phi_{0,k_2}(y) \,\mathrm{d}x
\,\mathrm{d}y,\\
\widehat\beta^{(c,e)}_{j;k_1,k_2}&\displaystyle\ninfty&\beta
^{(c,e)}_{j;k_1,k_2}:=\iint_{\R^d\times\R^d}h_c(x,y,\theta) \Psi
_{j;k_1,k_2}^{(e)}(x,y) \,\mathrm{d}x \,\mathrm{d}y
\end{eqnarray*}
on $(\XX_n^\theta)_n$. Hence,
$
\lim_{n\to\infty}n^2 \E[(\widehat U_{n,c}^{*(K,L)}-
U_{n,c}^{*(K,L)})^2|\X_n=x_n]= 0,
$
where the kernel of~$U_{n,c}^{*(K,L)}$ is obtained by substituting
$\widehat\alpha_{k_1,k_2}^{(c)}$ and $\widehat\beta
^{(c,e)}_{j;k_1,k_2}$ in the kernel of $\widehat U_{n,c}^{*(K,L)}$
through~$\alpha_{k_1,k_2}^{(c)}$ and $\beta^{(c,e)}_{j;k_1,k_2}$, respectively.

Thus, the next step is the application of the CLT of Neumann and Paparoditis \cite{NeuPa05}
to~$n U_{n,c}^{*(K,L)}$. For this purpose, we introduce
$Q_i^* :=\sum_{k=1}^{M(K,L)} t_k  q_k^*(X_i^*), t_1,\ldots,t_{M(K,L)}\in
\R$, where $q_k^*$ denotes the centered version (w.r.t.~$P_{X_1^*\vert\X
_n=x_n}$) of $\widetilde q_k$ and $(\widetilde q_k)_k$ is defined as in
the proof of Theorem~\ref{t.1}.
Obviously, given $X_1,\ldots,X_n$, the sequence $(Q_i^*)_{i}$ is
centered and has uniformly bounded second moments. Due to (A$1^*$)(i),
the Lindeberg condition is satisfied. In order to show that for
arbitrary $\varepsilon>0$ the inequalities $|\frac{1}{n}
\operatorname{var}(Q_1^*+\cdots+ Q_n^* |\X_n=x_n)-\sigma^2|<\varepsilon,  \forall
n\geq n_0(\varepsilon),$ hold true with $\sigma^2$ as in the proof of
Theorem~\ref{t.1}, the abbreviations
$\operatorname{var}^*(\cdot)= \operatorname{var}(\cdot|\X_n=x_n)$ and $\operatorname{cov}^*(\cdot)= \operatorname{cov}(\cdot|\X
_n=x_n)$ are used. Hence,
\begin{eqnarray*}
&& \biggl|\frac{1}{n}\operatorname{var}{}^*[Q_1^*+\cdots+Q_n^* ] - \sigma^2\biggr| \\
&&\quad \leq2\sum_{r=2}^\infty\min\biggl\{\frac{r-1}{n},1\biggr\}
|\operatorname{cov}{}^*(Q_1^*,Q_r^*)| +\Biggl|\operatorname{var}{}^*(Q_1^*)+2\sum_{r=2}^\infty
\operatorname{cov}{}^*(Q_1^*,Q_r^*) - \sigma^2\Biggr|\\
&&\quad \leq2\sum_{r=2}^\infty\min\biggl\{\frac{r-1}{n},1\biggr\}
|\operatorname{cov}{}^*(Q_1^*,Q_r^*)|
+2\Biggl|\sum_{r=2}^{R-1}[\operatorname{cov}{}^*(Q_1^*,Q_r^*)-\operatorname{cov}(Q_1,Q_r)]\Biggr|\\
&&\qquad{} +|\operatorname{var}{}^*(Q_1^*)-\operatorname{var}(Q_1)|+2\biggl|\sum_{r\geq R}
\operatorname{cov}{}^*(Q_1^*,Q_r^*)\biggr|
+2\biggl|\sum_{r\geq R} \operatorname{cov}(Q_1,Q_r)\biggr|.
\end{eqnarray*}
By (A1) and (A$1^*$), $R$ can be chosen such that $|\sum_{r\geq R}
\operatorname{cov}(Q_1,Q_r)| +|\sum_{r\geq R}\operatorname{cov}^*(Q_1^*, Q_r^*)|$
$\leq\varepsilon/4$. Moreover, (A$1^*$) implies that the first summand
can be bounded from above by~$\varepsilon/4$ as well if $n\geq
n_0(\varepsilon)$ for some $n_0(\varepsilon)\in\N$. According to the
convergence of the two-dimensional distributions and the uniform
boundedness of $(Q_k^*)_{k\in\Z}$, it is possible to pick
$n_0(\varepsilon)$ such that additionally the two remaining summands
are bounded by $\varepsilon/8$. For the\vadjust{\goodbreak} validity of the CLT of Neumann and Paparoditis \cite
{NeuPa05} in probability, it remains to verify
their inequality~(6.4). By Lipschitz continuity of $Q_{t_1}^*Q_{t_2}^*$
this holds with $\bar\theta_r=\operatorname{Lip}(Q_{t_1}^*Q_{t_2}^*)\bar\tau_r\leq
C\bar\tau_r$.
The application of the continuous mapping theorem results in
$n  U_{n,c}^{*(K,L)} \stackrel{d}{\longrightarrow} Z_c^{(K,L)}$, in
probability. Invoking the same arguments as in the proof of Theorem~\ref
{t.1}, this implies $n U_{n}^*\stackrel{d}{\longrightarrow} Z$, in probability.

In order to obtain the analogous result of convergence for $n V_n^*$,
we define \mbox{$\widetilde\XX_n^\theta\,{\subseteq}\,\XX_n^\theta, n\,{\in}\,\N,$}
such that
$|\E(|h(X_1^*,X_1^*,\widehat\theta_n)|\vert\X_n=x_n)-\E|h(X_1,X_1,\theta
)||\leq\eta_n, \forall x_n\in\widetilde\XX_n^\theta$. Here, the null
sequence $(\eta_n)_{n\in\N}$ is chosen in such a way that $P(\X_n\in
\widetilde\XX_n^\theta)\longrightarrow_{n\to\infty} 1.$
Now, additionally to our previous considerations,
\[
P\Biggl(\Biggl|\frac{1}{n}\sum_{i=1}^n h(X_i^*,X_i^*,\widehat\theta
_n)-\E h(X_1,X_1,\theta)\Biggr| > \varepsilon\Big|\X_n=x_n\Biggr)
\mathop{\longrightarrow} _{n\to\infty} 0
\]
has to be proved for arbitrary $\varepsilon>0$ and any sequence
$(x_n)_{n\in\N}$ with $x_n\in\widetilde\XX_n^\theta, n\in\N$.
According to the definition of the sets $(\widetilde\XX_n^\theta)_n$,
we get $\E(h(X_1^*,X_1^*,\widehat\theta_n)\vert\X_n=x_n
)\longrightarrow_{n\to\infty} \E h(X_1,X_1,\theta)$. Therefore, it
suffices to prove
\[
P\Biggl(\Biggl|\frac{1}{n}\sum_{k=1}^n \bigl[h(X_k^*,X_k^*,\widehat
\theta_n)-\E\bigl(h(X_1^*,X_1^*,\widehat\theta_n)|\X_n=x_n\bigr)\bigr]\Biggr| >
\frac{\varepsilon}{2} \Big|\X_n=x_n\Biggr)\ninfty0.
\]
This in turn is a consequence of Lemma~\ref{l.lln}
since under the assumptions of the theorem the sequence of functions
$(g_n)_{n\in\N}$ with $g^{(n)}(\cdot)=h(\cdot,\cdot,\widehat\theta_n)-\E
( h(X_1^*,X_1^*,\widehat\theta_n)\vert\X_n=x_n)$ is uniformly integrable
and satisfies the smoothness property presumed in Lemma~\ref{l.lln}.
Finally, bootstrap consistency follows from Lemma~\ref{l.6}.
\end{pf*}

\subsection{Proofs of auxiliary results}\label{SS52}

First, we derive a weak law of large numbers for smooth functions of
triangular arrays of $\tau$-dependent random variables.

\begin{lem}[(Weak law of large numbers)]\label{l.lln}
Let $(X_{n,k})_{k=1}^n, n\in \N,$ be a triangular scheme of
(row-wise) stationary, $\R^d$-valued, integrable random variables such
that $\lim_{K\to\infty} \sup_{n\in\N}P(\|X_{n,1}\|_{l_1}>K)=0.$
Suppose that the coefficients $\bar\tau_r:=\sup_{n>r}\tau_{r,n}$
satisfy $\bar\tau_r\longrightarrow_{r\to\infty}0$, where
\begin{eqnarray*}
\tau_{r,n}&:=&
\sup\{\tau(\sigma(X_{n,s_1},\ldots,X_{n,s_u}),(X_{n,t_1}^\prime
,X_{n,t_2}^\prime,X_{n,t_3}^\prime)^\prime)\vert u\in\N,\\
&&\phantom{\sup\{} 1\leq s_1\leq\cdots\leq s_u < s_u+r\leq t_1\leq t_2\leq t_3\leq n\}.
\end{eqnarray*}
Moreover, suppose that the functions $g^{(n)}\dvtx\R^d\to\R^p$ with $\E
g^{(n)}(X_{n,1})=0_p$ are uniformly Lipschitz continuous on any bounded
interval.
If additionally the sequence $(g^{(n)}(X_{n,1}))_{n\in\N}$ is uniformly
integrable, then
\[
\frac{1}{n}\sum_{k=1}^n g^{(n)}(X_{n,k})\stackrel{P}{\longrightarrow}
0_p.\vadjust{\goodbreak}
\]
\end{lem}

\begin{pf}
W.l.o.g.~let $p=1$. We prove that for arbitrary $\varepsilon, \eta>0$
there exists an $n_0$ such that for all $n>n_0$ the inequality
$P(|n^{-1}\sum_{k=1}^n g^{(n)}(X_{n,k})|>\varepsilon)\leq\eta$ holds.
To this end, a truncation argument is invoked. Let $w_K$ denote a
Lipschitz continuous, nonnegative function that is bounded from above
by one such that $w_K(x)=1$ for $x\in[-K,K]^d$ and $w_K(x)=0$ for
$x\notin[-K-1,K+1]^d$ with $K\in\R_+$. For a finite constant $M$, that
is specified later, define functions $ g_{M,K}^{(n)}\dvtx \R^d\to\R$ by
\[
g_{M,K}^{(n)}(x):=
\cases{
g^{(n)}(x) w_K(x)&\quad$\mbox{for }\bigl|g^{(n)}(x) w_K(x)\bigr|\leq
M,$\vspace*{1pt}\cr
-M&\quad$\mbox{for }g^{(n)}(x) w_K(x)< -M,$\vspace*{1pt}\cr
M&\quad$\mbox{for }g^{(n)}(x) w_K(x)> M$}\vspace*{-2pt}
\]
and $g_{M,K}^{(n,c)}$ by $g_{M,K}^{(n,c)}(x)= g_{M,K}^{(n)}(x)-\E
g^{(n)}_{M,K}(X_{n,1})$.
This allows for the estimation
\begin{eqnarray*}
P\Biggl(\Biggl|\frac{1}{n}\sum_{k=1}^n g^{(n)}(X_{n,k})\Biggr|>\varepsilon\Biggr)
&\leq& P\Biggl(\Biggl|\frac{1}{n}\sum_{k=1}^n g^{(n)}(X_{n,k})-
g_{M,K}^{(n)}(X_{n,k})\Biggr|>\frac{\varepsilon}{3}\Biggr)
\\[-3pt]
&&{} +P\biggl(\bigl|\E g_{M,K}^{(n)}(X_{n,1})\bigr|>\frac{\varepsilon}{3}
\biggr)+P\Biggl(\Biggl|\frac{1}{n}\sum_{k=1}^n g_{M,K}^{(n,c)}(X_{n,k})\Biggr|>\frac
{\varepsilon}{3}\Biggr).\vspace*{-2pt}
\end{eqnarray*}
According to Markov's inequality, the first summand on the r.h.s.~can
be bounded by
\[
\frac{3}{\varepsilon}\Bigl[\sup_{n\in\N}\E\bigl|g^{(n)}(X_{n,1})\bigr|\I
_{|g^{(n)}(X_{n,1})|>M}+M\sup_{n\in\N}P(\|X_{n,1}\|_{l_1}>K)\Bigr].\vspace*{-2pt}
\]
Since the functions $g^{(n)}, n\in\N,$ are centered, we additionally obtain
\begin{eqnarray*}
&&P\biggl(\bigl|\E g_{M,K}^{(n)}(X_{n,1})\bigr|>\frac{\varepsilon}{3}\biggr)\\[-3pt]
&&\quad\leq P\biggl(\sup_{n\in\N}\E\bigl|
g_{M,K}^{(n)}(X_{n,1})-g^{(n)}(X_{n,1})\bigr|>\frac{\varepsilon}{3}\biggr)\\[-3pt]
&&\quad \leq P\biggl(\sup_{n\in\N}\E\bigl|g^{(n)}(X_{n,1})\bigr|\I
_{|g^{(n)}(X_{n,1})|>M}+M\sup_{n\in\N}P(\|X_{n,1}\|_{l_1}>K)>\frac
{\varepsilon}{3}\biggr).\vspace*{-2pt}
\end{eqnarray*}
Therefore, by choosing $M$ and $K=K(M)$ sufficiently large,
we get
\[
P\Biggl(\Biggl|\frac{1}{n}\sum_{k=1}^n g^{(n)}(X_{n,k})-
g_{M,K}^{(n)}(X_{n,k})\Biggr|>\frac{\varepsilon}{3}\Biggr)
+P\biggl(\bigl|\E g_{M,K}^{(n)}(X_{n,1})\bigr|>\frac{\varepsilon}{3}
\biggr)\leq\frac{\eta}{2}.\vspace*{-2pt}
\]
Concerning the remaining term, Chebyshev's inequality leads to
\[
P\Biggl(\Biggl|\frac{1}{n}\sum_{k=1}^n g_{M,K}^{(n,c)}(X_{n,k})\Biggr|>\frac
{\varepsilon}{3}\Biggr)
\leq\frac{9M^2}{\varepsilon^2 n}+\frac{18}{\varepsilon^2 n^2}\sum
_{j<k}\E g_{M,K}^{(n,c)}(X_{n,j}) g_{M,K}^{(n,c)}(X_{n,k}).\vspace*{-2pt}
\]
Thus, it remains to derive an upper bound for $n^{-2}\sum_{j<k}|\E
g_{M,K}^{(n,c)}(X_{n,j}) g_{M,K}^{(n,c)}(X_{n,k})|$ that
vanishes\vadjust{\goodbreak}
asymptotically. For this purpose, we introduce a copy $\widetilde
X_{n,k}$ of $X_{n,k}$, that is independent of $ X_{n,j}$ and such that
$\E\|X_{n,k}-\widetilde X_{n,k}\|_{l_1}\leq\tau_{k-j,n}$. Due to their
construction, the functions $g_{M,K}^{(n,c)}$ are Lipschitz continuous
uniformly in $n$ and with a constant $C(M,K)$. This implies
\begin{eqnarray*}
\frac{1}{n^2}\sum_{j<k}\bigl|\E g_{M,K}^{(n,c)}(X_{n,j})
g_{M,K}^{(n,c)}(X_{n,k})\bigr|&\leq& \frac{2  M}{n^2}\sum_{j<k}\E\bigl|
g_{M,K}^{(n,c)}(X_{n,k})- g_{M,K}^{(n,c)}(\widetilde X_{n,k})\bigr|\\
&\leq&\frac{2 M C(M,K)}{n}\sum_{r=1}^n \bar\tau_r,
\end{eqnarray*}
where the remaining term converges to zero according to Cauchy's limit
theorem, cf. Knopp~\cite{Kn56}.
\end{pf}

In order to prove Lemma~\ref{l.1}, Lemma~\ref{l.5}, and Lemma~\ref
{l.4}, an approximation of terms of the structure
\[
Z_n :=\frac{1}{n^2}\mathop{\sum_{i,j,k,l=1}}_{i\neq j;k\neq l}^n \E
H(X_i,X_j)H(X_k,X_l)
\]
is required. Here, $H$ denotes a symmetric, degenerate kernel function.
Assuming that $(X_n)_{n\in\N}$ satisfies (A1), we obtain
\[
Z_n\leq\frac{8}{n^2}\sum_{i< j;k< l;i\leq k}^n|\E
H(X_i,X_j)H(X_k,X_l)|\\
\leq8\sup_{1\leq k<n}\E|H(X_1,X_{1+k})|^2+\frac{8}{n^2}\sum
_{r=1}^{n-1}\sum_{t=1}^4 Z_{n,r}^{(t)}
\]
with
\begin{eqnarray*}
Z_{n,r}^{(1)}&:=&\mathop{\sum_{1\leq i< j;k< l;j\leq l\leq n}}_{r:=\min\{
j,k\}-i\geq l-\max\{j,k\}}\bigl|\E H(X_i,X_j)H(X_k,X_l)-\E
H\bigl(X_i,\widetilde X_j^{(r)}\bigr)H\bigl(\widetilde X_k^{(r)}, \widetilde
X_l^{(r)}\bigr)\bigr|,\\
Z_{n,r}^{(2)}&:=&\mathop{\sum_{1\leq i< j;i\leq k;k< l\leq n}}_{r:=l-\max
\{j,k\}>\min\{j,k\}-i}\bigl|\E H(X_i,X_j)H(X_k,X_l)-\E
H(X_i,X_j)H\bigl(X_k,\widetilde X_l^{(r)}\bigr)\bigr|,\\
Z_{n,r}^{(3)}&:=&\mathop{\sum_{1\leq i\leq k< l< j\leq n}}_{r:=k-i\geq
j-l}\bigl|\E H(X_i,X_j)H(X_k,X_l)-\E H\bigl(X_i,\widetilde
X_j^{(r)}\bigr)H\bigl(\widetilde X_k^{(r)},\widetilde X_l^{(r)}\bigr)\bigr|,\\
Z_{n,r}^{(4)}&:=&\mathop{\sum_{1\leq i\leq k< l< j\leq n}}_{
r:=j-l>k-i}\bigl|\E H(X_i,X_j)H(X_k,X_l)-\E H\bigl(X_i,\widetilde
X_j^{(r)}\bigr)H(X_k, X_l)\bigr|.
\end{eqnarray*}
Here, in every summand of $Z_{n,r}^{(1)}$ and $Z_{n,r}^{(3)}$ the
vector $(\widetilde X_j^{(r)\prime},\widetilde X_k^{(r)\prime
},\widetilde X_l^{(r)\prime})^\prime$ is chosen such that it is
independent of the random variable $X_i$, $(\widetilde X_j^{(r)\prime
},\widetilde X_k^{(r)\prime},\widetilde X_l^{(r)\prime})^\prime\stackrel
{d}{=}( X_j^\prime, X_k^\prime, X_l^\prime)^\prime$, and (\ref{eq.a1})
holds. Within $Z_{n,r}^{(2)}$ (resp., $Z_{n,r}^{(4)})$, the
random variable $\widetilde X_l^{(r)}$ (resp., $\widetilde
X_j^{(r)}$) is chosen to be independent of the vector $( X_i^\prime,
X_j^\prime, X_k^\prime)^\prime$ (resp., $( X_i^\prime, X_k^\prime
, X_l^\prime)^\prime$) such that $ \widetilde X_l^{(r)}\stackrel{d}{=}
X_l$ (resp., $ \widetilde X_j^{(r)}\stackrel{d}{=} X_j$) and
(\ref{eq.a1}) holds. This may possibly require an enlargement of the
underlying probability space.
Moreover, note that the subtrahends of these expressions vanish due to
the degeneracy of $H$ and that the number of summands of
$Z_{n,r}^{(t)}, t=1,\ldots,4,$ is bounded by $(r+1)n^2$.
For sake of notational simplicity, the upper index $r$ is omitted in
the sequel.

\begin{pf*}{Proof of Lemma~\textup{\protect\ref{l.1}}}
For $c>0$, we define $c_h:=\max_{x,y\in[-c,c]^d} |h(x,y)|$,
\[
\widetilde h^{(c)}(x,y) :=
\cases{
h(x,y) &\quad$\mbox{for } |h(x,y)|\leq c_h$,\vspace*{2pt}\cr
-c_h &\quad$\mbox{for } h(x,y) <-c_h,$\vspace*{2pt}\cr
c_h &\quad$\mbox{for } h(x,y) >c_h$}
\]
and its degenerate version
\begin{eqnarray*}
h_c(x,y) &:= &\widetilde h^{(c)}(x,y) - \int_{\R^d}\widetilde
h^{(c)}(x,y) P_X(\mathrm{d}x)
- \int_{\R^d}\widetilde h^{(c)}(x,y) P_X(\mathrm{d}y)\\
&&{}+ \iint_{\R^d\times\R^d} \widetilde h^{(c)}(x,y) P_X(\mathrm{d}x) P_X(\mathrm{d}y).
\end{eqnarray*}
The approximation error $n^2 \E(U_n-U_{n,c})^2$ can be reformulated in
terms of $Z_n$ with kernel $H=H^{(c)}:=h-h^{(c)}$. Hence, it remains to
verify that $\sup_{k\in\N}\E|H^{(c)}(H_1,X_{1+k})|^2$ and $\sup_{n\in\N
}n^{-2}\sum_{r=1}^{n-1}\sum_{t=1}^4 Z_{n,r}^{(t)}$
tend to zero as $c\to\infty.$
First, we consider $\sup_{n\in\N}n^{-2}\sum_{r=1}^{n-1}$
$Z_{n,r}^{(1)}$, the remaining quantities can be treated similarly. The
summands of $Z_{n,r}^{(1)}$ are bounded as follows:
%
\begin{eqnarray}\label{eq.unb1}
&& \bigl| \E H^{(c)}(X_i,X_j)H^{(c)}(X_k, X_l)- \E H^{(c)}(X_i,\widetilde
X_j)H^{(c)}(\widetilde X_{k},\widetilde X_{l})\bigr|\nonumber\\
&&\quad\leq\E\bigl|H^{(c)}(X_k,X_l)
\bigl[H^{(c)}(X_i,X_j)-H^{(c)}(X_i,\widetilde X_{j})\bigr]\I_{(X_k^\prime
,X_l^\prime)^\prime\in[-c,c]^{2d}}\bigr|\nonumber\\
&&\qquad{} +\E\bigl|H^{(c)}(X_k,X_l)
\bigl[H^{(c)}(X_i,X_j)-H^{(c)}(X_i,\widetilde X_{j})\bigr]\I_{(X_k^\prime,
X_l^\prime)^\prime\notin[-c,c]^{2d}}\bigr|
\nonumber
\\[-8pt]
\\[-8pt]
\nonumber
&&\qquad{} + \E\bigl|H^{(c)}(X_i,\widetilde X_j)
\bigl[H^{(c)}(X_k,X_l)-H^{(c)}(\widetilde X_{k},\widetilde X_{l})\bigr]\I
_{(X_i^\prime,\widetilde X_j^\prime)^\prime\in[-c,c]^{2d}}\bigr|\\
&&\qquad{} + \E\bigl|H^{(c)}(X_i,\widetilde X_j)
\bigl[H^{(c)}(X_k,X_l)-H^{(c)}(\widetilde X_{k},\widetilde X_{l})\bigr]\I
_{(X_i^\prime,\widetilde X_j^\prime)^\prime\notin[-c,c]^{2d}}\bigr|\nonumber\\
&&\quad= E_1+E_2+E_3+E_4.\nonumber
\end{eqnarray}

The functions $H^{(c)}$ are obviously Lipschitz continuous uniformly in $c$.
Therefore, an iterative application of H\"older's inequality to $E_2$
yields
%
\begin{eqnarray}\label{eq.zerl2}
E_2&\leq& \bigl(\E\bigl|H^{(c)}(X_i,X_j)-H^{(c)}(X_i,\widetilde
X_{j})\bigr|\bigr)^\delta\nonumber\\
&&{}\times\bigl(\E\bigl|H^{(c)}(X_k,X_l)\bigr|^{1/(1-\delta
)}\bigl|H^{(c)}(X_i,X_j)-H^{(c)}(X_i,\widetilde X_{j})\bigr|\I_{(X_k^\prime,
X_l^\prime)^\prime\notin[-c,c]^{2d}}\bigr)^{1-\delta}\qquad\nonumber\\
&\leq&  C\tau_r^\delta\bigl\{\bigl(\E\bigl|H^{(c)}(X_k,X_l)\bigr|^{(2-\delta
)/(1-\delta)}\I_{(X_k^\prime, X_l^\prime)^\prime\notin
[-c,c]^{2d}}\bigr)^{1/(2-\delta)}\\
&&{}\times\bigl(\E\bigl|H^{(c)}(X_i,X_j)\bigr|^{(2-\delta)/(1-\delta)}+\E
\bigl|H^{(c)}(X_i,\widetilde X_j)\bigr|^{(2-\delta)/(1-\delta)}\bigr)^{(1-\delta
)/(2-\delta)}\bigr\}^{1-\delta}\nonumber\\
&\leq& C\tau_r^\delta\bigl(\E\bigl|H^{(c)}(X_k,X_l)\bigr|^{(2-\delta)/(1-\delta
)}\I_{(X_k^\prime, X_l^\prime)^\prime\notin[-c,c]^{2d}}
\bigr)^{(1-\delta)/(2-\delta)}.\nonumber
\end{eqnarray}
As $\sup_{k\in\N}\E| h(X_1,X_{1+k})|^\nu<\infty$ for $\nu>(2-\delta
)/(1-\delta)$, we obtain $E_2\leq\tau_r^\delta \varepsilon_1(c)$ with
$\varepsilon_1(c)\longrightarrow_{c\to\infty}0$ after employing H\"
older's inequality once again.
Analogous calculations yield $E_4\leq\tau_r^\delta \varepsilon_2(c)$
with $\varepsilon_2(c)\longrightarrow_{c\to\infty}0$. Likewise, the
approximation methods for $E_1$ and $E_3$ are equal. Therefore, only
$E_1$ is considered:
\begin{eqnarray*}
E_1&\leq& \E\biggl|\int_{\R^d}\widetilde h^{(c)}(X_k,y)P_X(\mathrm{d}y)
\bigl[H^{(c)}(X_i,X_j)-H^{(c)}(X_i,\widetilde X_{j})\bigr]\I_{X_k\in
[-c,c]^d}\biggr|\\
&&{}+\E\biggl|\int_{\R^d}\widetilde h^{(c)}(y,X_l)P_X(\mathrm{d}y)
\bigl[H^{(c)}(X_i,X_j)-H^{(c)}(X_i,\widetilde X_{j})\bigr]\I_{X_l\in
[-c,c]^d}\biggr|\\
&&{}+\E\biggl|\iint_{\R^d\times\R^d}\widetilde
h^{(c)}(x,y)P_X(\mathrm{d}x)P_X(\mathrm{d}y)\bigl[H^{(c)}(X_i,X_j)-H^{(c)}(X_i,\widetilde
X_{j})\bigr]\biggr|\\
&=& E_{1,1}+E_{1,2}+E_{1,3}.
\end{eqnarray*}
Analogous to ~(\ref{eq.zerl2}), we obtain
\begin{eqnarray*}
E_{1,1}&\leq& C \tau_r^\delta \biggl\{\biggl(\E\biggl|\int_{\R^d}
h(X_k,y)-\widetilde h^{(c)}(X_k,y)P_X(\mathrm{d}y)\biggr|^{(2-\delta)/(1-\delta)}\I
_{X_k\in[-c,c]^d}\biggr)^{1/(2-\delta)}\\
&&{}\times\Bigl[\sup_{k\in\N}\E\bigl|H^{(c)}(X_1,X_{1+k})\bigr|^{(2-\delta)/(1-\delta
)}+\E\bigl|H^{(c)}(X_i,\widetilde X_j)\bigr|^{(2-\delta)/(1-\delta)}
\Bigr]^{(1-\delta)/(2-\delta)}\biggr\}^{1-\delta}\\
&\leq& C \tau_r^\delta\biggl (\int_{\R^d}\int_{\R^d}
\bigl|h(x,y)-\widetilde h^{(c)}(x,y)\bigr|^{(2-\delta)/(1-\delta)}\\
&&\hspace*{54pt}{}\times P_X(\mathrm{d}y)\I_{x\in
[-c,c]^d} P_X(\mathrm{d}x)\biggr)^{(1-\delta)/(2-\delta)}\\
&\leq&  \tau_r^\delta \varepsilon_3(c)
\end{eqnarray*}
with $\varepsilon_3(c)\longrightarrow_{c\to\infty}0$. The estimation of
$E_{1,2}$ coincides with the previous one.
The expression $E_{1,3}$ can be bounded as follows:
\begin{eqnarray*}
E_{1,3}&\leq& C \tau_r\iint_{\R^d\times\R^d}\bigl|h(x,y)-\widetilde
h^{(c)}(x,y)\bigr|P_X(\mathrm{d}x)P_X(\mathrm{d}y)\\
&\leq& C \tau_r \iint_{\R^d\times\R^d}|h(x,y)|\I_{(x^\prime,y^\prime
)^\prime\notin[-c,c]^{2d}}P_X(\mathrm{d}x)P_X(\mathrm{d}y)\\
&\leq&\tau_r \varepsilon_4(c)
\end{eqnarray*}
with $\varepsilon_4(c)\longrightarrow_{c\to\infty}0.$ To sum up,
we have $E_1+E_2+E_3+E_4\leq\varepsilon_5(c) \tau_r^\delta$, where
$\varepsilon_5(c)\longrightarrow_{c\to\infty}0$ uniformly in $n$. This
leads to
\[
\lim_{c\to\infty}\sup_{n\in\N}\frac{1}{n^2}\sum_{r=1}^{n-1}Z_{n,r}^{(1)}
\leq\lim_{c\to\infty}\sup_{n\in\N} \frac{1}{n^2}\sum_{r=1}^{n-1}(r+1)
n^2 \tau_r^\delta \varepsilon_5(c)
=0.
\]
It remains to examine
\begin{eqnarray*}
\sup_{k\in\N}\E\bigl[H^{(c)}(X_1,X_{1+k})\bigr]^2 &\leq& C\Bigl( \sup_{k\in\N} \E
\bigl[h(X_1,X_{1+k})-\widetilde h^{(c)}(X_1,X_{1+k})\bigr]^2\\[-3pt]
&&\hphantom{C\Bigl(}{}+ \E\bigl[h(X_1,\widetilde X_1)-\widetilde h^{(c)}(X_1,\widetilde
X_1)\bigr]^2\Bigr).
\end{eqnarray*}
Here, $\widetilde X_1$ denotes an independent copy of $X_1$. Similar
arguments as before yield $\lim_{c\to\infty}\sup_{k\in\N}\E
[H^{(c)}(X_1,X_{1+k})]^2=0.$
\end{pf*}

The characteristics stated in the following two lemmas will be
essential for a wavelet approximation of the kernel function~$h$.

\begin{lem}\label{l.2}
Given a Lipschitz continuous function $g\dvtx\R^d\to\R$, define a wavelet
series approximation~$g_j$ by $g_j(x):=\sum_{k\in\Z^d}\alpha_{j,k}\Phi
_{j,k}(x),j\in\Z$, where $\alpha_{j,k}=\int_{\R^d}g(x) \Phi_{j,k}(x)
\,\mathrm{d}x$. Then $g_j$ is Lipschitz continuous with a constant that is
independent of $j$.
\end{lem}

\begin{pf}
In order to establish Lipschitz continuity, the function~$g_{j}$ is
decomposed into two parts
\begin{eqnarray*}
g_{j}(x)
&= &\sum_{ k \in\Z^d} \biggl[ \int_{\R^d}\Phi_{j,k}(u) g(x)\, \mathrm{d}u\biggr]
\Phi_{j,k}(x)
+ \sum_{ k \in\Z^d} \biggl[ \int_{\R^d} \Phi_{j,k}(u) [g(u) - g(x)]
\,\mathrm{d}u\biggr]
\Phi_{j,k}(x) \\
&=& H_1(x) + H_2(x).
\end{eqnarray*}
According to the above choice of the scale function (with
characteristics (1)--(3) of Section~\ref{SS22}), the prerequisites
of Corollary~8.1 of H{\"a}rdle \textit{et al.} \cite{Hetal98} are fulfilled for $N=1$. This
implies that $\int_{-\infty}^\infty\sum_{l\in\Z} \phi(y-l)\phi(z-l)\,\mathrm{d}z
=1, \forall  y\in\R$. Based on this result, we obtain
\[
\sum_{ k \in\Z^d} \int_{\R^d}\Phi_{j,k}(u)\Phi_{j,k}(x)  \,\mathrm{d}u
= 2^{jd} \prod_{i=1}^d\int_{-\infty}^\infty\sum_{l \in\Z}\phi
(2^{j}u_i-l)\phi(2^{j}x_i-l)\,\mathrm{d}u_i=1\qquad \forall  x\in\R^d,
\]
by applying an appropriate variable substitution. To this end, note
that for every fixed $x$, the number of non-vanishing summands can be
bounded by a finite constant uniformly in~$j$ because of the finite
support of $\phi$. Therefore, the order of summation and integration is
interchangeable. Hence, $H_1=g$ which in turn immediately implies the
desired continuity property for~$H_1$.

In order to investigate $H_2$, we define a sequence of functions
$(\kappa_{k})_{k\in\Z}$ by
\[
\kappa_{k}(x) = \int_{\R^d}\Phi_{j,k}(u) [g(u) - g(x)]
\,\mathrm{d}u.\vadjust{\goodbreak}
\]
These functions are Lipschitz continuous with a constant decreasing in $j$:
%
\begin{equation}\label{eq.1}
| \kappa_{k}(x) - \kappa_{k}(\bar x) |
\leq \operatorname{Lip}(g)  \mathrm{O}(2^{-jd/2})  \|x-\bar x\|_{l_1}.
\end{equation}
Moreover, boundedness and Lipschitz continuity of~$\phi$ yield
%
\begin{equation}\label{eq.2}
\| \Phi_{j,k} \|_\infty= \mathrm{O}(2^{jd/2})\quad \mbox{and}\quad
|\Phi_{j,k}(x) - \Phi_{j,k}(\bar x) |
= \mathrm{O}\bigl(2^{j(d/2+1)}\bigr) \|x-\bar x\|_{l_1}.
\end{equation}
Thus,
\begin{eqnarray*}
|H_2(x)-H_2(\bar x)|
&\leq &\sum_{k\in\Z^d} |\Phi_{j,k}(x)| |\kappa_{k}(x) - \kappa
_{k}(\bar x) |\\
&&{}+ \sum_{ k \in\Z^d} |\kappa_{k}(\bar x)||\Phi
_{j,k}(x) - \Phi_{j,k}(\bar x) |\\
&\leq&C \|x-\bar x\|_{l_1} + \sum_{k\in\Z^d} |\kappa_{k}(\bar x)|
|\Phi_{j,k}(x) - \Phi_{j,k}(\bar x) |.
\end{eqnarray*}
Now, it has to be distinguished whether or not $\bar x\in\operatorname{supp}(\Phi_{j,k})$
in order to approximate the second summand. (Here,
$\operatorname{supp}$
denotes the support of a function.) In the first case, it is
helpful to illuminate $|\kappa_{k}(\bar x)|= |\int_{\R^d} \Phi_{j,k}(u)
[g(u) - g(\bar x)] \,\mathrm{d}u|$. The integrand is non-trivial only if $u\in
\operatorname{supp } (\Phi_{j,k})$. In these situations, $|g(u) - g(\bar
x)|=\mathrm{O}(2^{-j})$ by Lipschitz continuity. Consequently, we get
\[
|\kappa_{k}(\bar x)|\leq \mathrm{O}(2^{-j})\int_{\R^d}|\Phi_{j,k}(u)|\,\mathrm{d}u= \mathrm{O}\bigl(2^{-j(d/2+1)}\bigr)
\]
which leads to
\[
\sum_{k\in\Z^d} |\kappa_{k}(\bar x)||\Phi_{j,k}(x) - \Phi
_{j,k}(\bar x) |\leq C \|x-\bar x\|_{l_1}
\]
as the number of nonvanishing summands is finite, independently of the
values of $x$ and~$\bar x$. Therefore, Lipschitz continuity of $H_2$ is
obtained as long as $\bar x\in\operatorname{supp }(\Phi_{j,k})$.

In the opposite case, we only have to consider the situation of $x\in
\operatorname{supp }(\Phi_{j,k})$ since the setting~$\bar x$, $x\notin\operatorname
{supp }(\Phi_{j,k})$ is trivial. With the aid of (\ref{eq.1}) and (\ref
{eq.2}), the first term of the r.h.s.~of
%
\begin{equation}\label{eq.l2.1}
|\kappa_{k}(\bar x)[\Phi_{j,k}(x) - \Phi_{j,k}(\bar x)
] |\leq|\kappa_{k}(\bar x) - \kappa_{k}(x) | |\Phi
_{j,k}(x)|+ |\kappa_{k}(x) | |\Phi_{j,k}(x) - \Phi_{j,k}(\bar x)
|
\end{equation}
can be estimated from above by $C\|x-\bar x\|_{l_1} $. The
investigation of the second summand is identical to the analysis of the
case $\bar x\in\operatorname{supp }(\Phi_{j,k})$.

Finally, we obtain $|H_2(x)-H_2(\bar x)| \leq C \|x-\bar x\|_{l_1}$,
where $ C < \infty$ is a constant that is independent of $j$. This
yields the assertion of the lemma.
\end{pf}

\begin{lem}\label{l.3}
Let $g\dvtx \R^d\to\R$ be a function that is continuous on some interval $(-c,c)^d$.
For arbitrary $b\in(0,c)$ and $K\in\N$ there exists a $J(K,b,c)\in\N$
such that for $g$ and its approximation~$ g_J$ given by $ g_J(x) = \sum
_{ k \in\Z^d} \alpha_{J,k}\Phi_{J,k}(x)$ it holds
\[
\max_{x\in[-b,b]^d} |g(x)  -  {g}_J(x)|
 \leq  1/K\qquad \forall  J\geq J(K,b,c) .\vspace*{-2pt}
\]
\end{lem}

\begin{pf}
Given $b\in(0,c)$, we define $\bar g^{(b,c)}(x):=g(x) w_{b,c}(x)$,
where $w_{b,c}$ is a Lipschitz continuous and nonnegative weight
function with compact support $S_w\subset(-c,c)^d$. Moreover,~$w_{b,c}$
is assumed to be bounded from above by~1 and $w_{b,c}(x) :=1$ for $x
\in(-b-\delta,\allowbreak b+\delta)^d$ for some $\delta>0$ with $b+\delta< c$.
Additionally, we set
$\alpha^{(b,c)}_{J,k}:=\int_{\R^d} \bar g^{(b,c)}(u)\Phi_{J,k}(u) \,\mathrm{d}u$. Hence,
\begin{eqnarray*}
&&\max_{ x\in[-b,b]^d}|g(x)  -  g_J(x)|\\[-2pt]
&&\qquad\leq \max_{ x\in[-b,b]^d} \biggl|\bar g^{(b,c)} (x)
- \sum_{k \in\Z^d} \alpha^{(b,c)}_{J, k} \Phi_{J,k}(x) \biggr|+ \max_{
x\in[-b,b]^d}\biggl |\sum_{k \in\Z^d} \alpha^{(b,c)}_{J, k} \Phi_{J,k}(x)
 -  g_J(x)\biggr|\\[-2pt]
&&\qquad= \max_{ x\in[-b,b]^d} A^{(J)}(x)+ \max_{ x\in[-b,b]^d}B^{(J)}(x).
\end{eqnarray*}
Since $\bar g^{(b,c)}\in C_0(\R^d)$, Theorem~8.4 of Wojtaszczyk \cite{Woj97}
implies that there exists a~$J_0(K,b,c)\in\N$ such that $\max_{ x\in
[-b,b]^d} A^{(J)}(x)\leq1/K$ for all $J\geq J_0(K,b,c)$.
Moreover, the introduction of the finite set of indices
\[
\bar Z(J):= \{k\in\Z^d  |  \Phi_{J,k}(x)\neq0 \mbox{ for some
} x\in[-b,b]^d\}
\]
leads to
\[
\max_{x\in[-b,b]^d}B^{(J)}(x) =
\max_{x\in[-b,b]^d}\biggl|
\sum_{k\in\bar Z(J)} \bigl( \alpha_{J,k} -\alpha_{J,k}^{(b,c)}\bigr) \Phi
_{J,k}(x)\biggr|.
\]
This term is equal to zero for all $J\geq J(K,b,c)$ and some
$J(K,b,c)\geq J_0(K,b,c)$ since the definition of $\bar g^{(b,c)}$
implies $\alpha_{J,k} = \alpha_{J,k}^{(b,c)},~ \forall k\in\bar Z,$
for all sufficiently large~$J$.
\end{pf}

\begin{pf*}{Proof of Lemma~\textup{\protect\ref{l.5}}}
The assertion of the lemma is verified in two steps.
First, the bounded kernel $h_c$, constructed in the proof of Lemma~\ref
{l.1}, is approximated by $\widetilde h^{(K)}_c$ which is defined by
$\widetilde h^{(K)}_c(x,y)=\sum_{k_1,k_2\in\mathbb{Z}^d}\alpha
_{J(K);k_1,k_2}^{(c)}\Phi_{J(K),k_1}(x)\Phi_{J(K),k_2}(y)$ with $\alpha
_{J(K);k_1,k_2}^{(c)}=\iint_{\R^d\times\R^d}h_c(x, y)\Phi
_{J(K),k_1}(x)\Phi_{J(K),k_2}(y) \,\mathrm{d}x \,\mathrm{d}y$. Here, the indices
$(J(K))_{K\in\N}$ with $J(K)$ $\longrightarrow_{K\to\infty} \infty$ are
chosen such that the assertion of Lemma~\ref{l.3} holds true for
$b=b(K)\in\R$ with $P(X_1\notin[-b,b]^d)\leq K^{-1}$ and $c=2b$.
Since the function $\widetilde h^{(K)}_c$ is not degenerate in general,
we introduce its degenerate counterpart
\begin{eqnarray*}
h^{(K)}_c(x,y)&= &\widetilde h^{(K)}_c(x,y)-\int_{\R^d}\widetilde
h^{(K)}_c(x,y)P_X(\mathrm{d}x)-\int_{\R^d}\widetilde h^{(K)}_c(x,y)P_X(\mathrm{d}y)\\[-2pt]
&&{}+\iint_{\R^d\times\R^d}\widetilde h^{(K)}_c(x,y)P_X(\mathrm{d}x)P_X(\mathrm{d}y)
\end{eqnarray*}
and denote the corresponding $U$-statistic by $U_{n,c}^{(K)}$.\vadjust{\goodbreak}

Now, the structure of the proof is as follows. First, we prove
%
\begin{equation}\label{eq.l5.1}
\sup_{n\in\N}n^2 \E\bigl(U_{n,c}-U_{n,c}^{(K)}\bigr)^2 \mathop
{\longrightarrow} _{K\to\infty} 0.
\end{equation}
In a second step, it remains to show that for every fixed $K$
%
\begin{equation}\label{eq.l5.2}
\sup_{n\in\N}n^2 \E\bigl(U_{n,c}^{(K)}-U_{n,c}^{(K,L)}\bigr)^2 \mathop
{\longrightarrow} _{L\to\infty} 0.
\end{equation}
In order to verify (\ref{eq.l5.1}), we rewrite $n^2 \E
(U_{n,c}-U_{n,c}^{(K)})^2$ in terms of $Z_n$ with kernel function
$H:=H^{(K)}=h_c-h^{(K)}_c$. Hence, it remains to verify that $\sup_{n\in
\N}n^{-2}\sum_{r=1}^{n-1}\sum_{t=1}^4 Z_{n,r}^{(t)}$
and $\sup_{k\in\N}\E|H^{(K)}(H_1,X_{1+k})|^2$ tend to zero as $K\to
\infty.$
Exemplarily, we investigate $\sup_{n\in\N}n^{-2}\sum
_{r=1}^{n-1}Z_{n,r}^{(1)}$. The summands of $Z_{n,r}^{(1)}$ can be
bounded as follows:
\begin{eqnarray*}
&&\bigr| \E H^{(K)}(X_i,X_j)H^{(K)}(X_k, X_l)-H^{(K)}(X_i,\widetilde
X_j)H^{(K)}(\widetilde X_k,\widetilde X_l)\bigr|\\
&&\quad\leq \E\bigl|H^{(K)}(X_k,X_l)
\bigl[H^{(K)}(X_i,X_j)-H^{(K)}(X_i,\widetilde X_{j})\bigr]\bigr|\\
&&\qquad{} + \E\bigl|H^{(K)}(X_i,\widetilde X_j)
\bigl[H^{(K)}(X_k,X_l)-H^{(K)}(\widetilde X_{k},\widetilde X_{l})
\bigr]\bigr|.
\end{eqnarray*}
Since further approximations are similar for both summands, we
concentrate on the first one. Note that boundedness of $h_c$ implies
uniform boundedness of $(H^{(K)})_K$ due to the compact support of the
function $\phi$. Moreover, the constant $\operatorname{Lip}(H^{(K)})$ does not depend
on $K$ in consequence of Lemma~\ref{l.2}. Therefore, the application of
H\"older's inequality leads to
\[
\E\bigl|H^{(K)}(X_k,X_l)\bigl[H^{(K)}(X_i,X_j)-H^{(K)}(X_i,\widetilde
X_{j})\bigr]\bigr|
\leq C\tau_r^\delta\bigl[\E\bigl|H^{(K)}(X_k,X_l)\bigr|^{1/(1-\delta)}
\bigr]^{1-\delta}.
\]
The construction of the sequence $(b(K))_K$ above allows for the
following estimation:
\begin{eqnarray*}
&&\E\bigl|H^{(K)}(X_k,X_l)\bigr|^{1/(1-\delta)}\\
&&\quad= \E\bigl|H^{(K)}(X_k,X_l)\bigr|^{1/(1-\delta)}\I_{X_k, X_l\in[-b(K),
b(K)]^d}+\mathrm{O}\bigl(P\bigl(X_1\notin[-b(K),b(K)]^d\bigr)\bigr)\\
&&\quad\leq\sup_{x,y\in\bigl[-b(K),b(K)\bigr]^d} \bigl|H^{(K)}(x,y)\bigr|^{1/(1-\delta)}+\frac{C}{K}.
\end{eqnarray*}
According to Lemma~\ref{l.3} and the above choice of the sequence
$(b(K))_K$, we obtain
\begin{eqnarray*}
&&\sup_{x,y\in[-b(K),b(K)]^d}\bigl|H^{(K)}(x,y)\bigr|\\
&&\quad\leq\frac{1}{K} + 2\sup_{{x,y\in[-b(K),b(K)]^d}} \E\bigl|h_c(x,
X_1)-\widetilde{h}^{(K)}_c(x,X_1)\bigr|\\
&&\qquad{} + \biggl| \iint_{\R^d\times\R^d}h_c(x,y)- \widetilde
{h}^{(K)}_c(x,y)   P_X(\mathrm{d}x)   P_X(\mathrm{d}y)\biggr|\\
&&\quad\leq\frac{4}{K}+ 2\sup_{{x\in[-b(K),b(K)]^d}} \E\bigl|h_c(x,
X_1)-\widetilde{h}^{(K)}_c(x,X_1)\bigr|\I_{X_1\notin[-b(K),b(K)]^d}\\
&&\qquad{} + 2 \int_{\R^d}\int_{\R^d\setminus[-b(K),b(K)]^d}
\bigl|h_c(x,y)-\widetilde{h}^{(K)}_c (x,y)\bigr|   P_X(\mathrm{d}x)   P_X(\mathrm{d}y)\\
&&\quad\leq\frac{C}{K}.
\end{eqnarray*}
Consequently,
\[
\bigl| \E H^{(K)}(X_i,X_j)H^{(K)}(X_k, X_l)-\E H^{(K)}(X_i,\widetilde
X_j)H^{(K)}(\widetilde X_k,\widetilde X_l)\bigr|\leq C  \varepsilon
_K  \tau^{\delta}_r
\]
for some null sequence $(\varepsilon_K)_K$. This implies that $\sup
_{n\in\N} n^{-2}\sum_{r=1}^n Z_{n,r}^{(1)}$ tends to zero as $K$ increases.
Furthermore, one obtains $\sup_{k\in\N}\E
[H^{(K)}(X_1,X_{1+k})]^2=\mathrm{O}(K^{-1})$ similarly to the
consideration of $\E|H^{(K)}(X_k,X_l)|^{1/(1-\delta)}$ above.
Thus, we get $\sup_{n}n^2 \E(U_{n,c}-U_{n,c}^{(K)})^2 \longrightarrow
_{K\to\infty} 0.$

The main goal of the previous step was the multiplicative separation of
the random variables which are cumulated in $h_c$.
The aim of the second step is the approximation of~$h_c^{(K)}$, whose
representation is given by an infinite sum, by a function consisting of
only finitely many summands.
Similar to the foregoing part of the proof the approximation error
$n^2 \E(U_{n,c}^{(K)}-U_{n,c}^{(K,L)})^2$ is reformulated in terms of
$Z_n$ with kernel $H:=H^{(L)}= h^{(K)}_c-h^{(K,L)}_c$. As before, we
exemplarily take $n^{-2}\sum_{r=1}^{n-1}Z_{n,r}^{(1)} $ and $\sup_{k\in
\N}\E|H^{(L)}(X_1,X_{1+k})|^2$ into further consideration.
Concerning the summands of $Z_{n,r}^{(1)}$, we obtain
\begin{eqnarray*}
&&\bigl| \E H^{(L)}(X_i,X_j)H^{(L)}(X_k, X_l)-\E H^{(L)}(X_i,\widetilde
X_j)H^{(L)}(\widetilde X_k,\widetilde X_l)\bigr|\\
&&\quad\leq\E\bigl|H^{(L)}(X_k,X_l)
\bigl[H^{(L)}(X_i,X_j)-H^{(L)}(X_i,\widetilde X_{j})\bigr]\I_{(X_k^\prime
,X_l^\prime)^\prime\in[-B,B]^{2d}}\bigr|\\
&&\qquad{} +\E\bigl|H^{(L)}(X_k,X_l)
\bigl[H^{(L)}(X_i,X_j)-H^{(L)}(X_i,\widetilde X_{j})\bigr]\I_{(X_k^\prime
,X_l^\prime)^\prime\notin[-B,B]^{2d}}\bigr|\\
&&\qquad{} + \E\bigl|H^{(L)}(X_i,\widetilde X_j)
\bigl[H^{(L)}(X_k,X_l)-H^{(L)}(\widetilde X_{k},\widetilde X_{l})\bigr]\I
_{(X_i^\prime,\widetilde X_j^\prime)^\prime\in[-B,B]^{2d}}\bigr|\\
&&\qquad{} + \E\bigl|H^{(L)}(X_i,\widetilde X_j)
\bigl[H^{(L)}(X_k,X_l)-H^{(L)}(\widetilde X_{k},\widetilde X_{l})\bigr]\I
_{(X_i^\prime,\widetilde X_j^\prime)^\prime\notin[-B,B]^{2d}}\bigr|\\
&&\quad=E_1+E_2+E_3+E_4
\end{eqnarray*}
for arbitrary $B>0$. Obviously, it suffices to take the first two
summands into further considerations. The both remaining terms can be
treated similarly. First, note that $(H^{(L)})_L$ is uniformly bounded.
Since $\phi$ and $\psi$ have compact support, the number of overlapping
functions within $(\Phi_{0,k})_{k\in\{-L,\ldots,L\}^d}$ and $(\Psi
_{j,k}^{(e)})_{k\in\{-L,\ldots,L\}^d, 0\leq j<J(K),e\in E}$ can be
bounded by a constant that is independent of $L$. By Lipschitz
continuity of $\phi$ and $\psi$, this leads to uniform Lipschitz
continuity of $(h^{(K,L)}_c)_{L\in\N}$.
Due to the reformulation
\[
\widetilde h^{(K)}_c(x,y)
= \sum_{ k_1,k_2\in\Z^d } \alpha^{(c)}_{k_1,k_2}\Phi_{0,k_1}(x)\Phi
_{0,k_2}(y)+ \sum_{j=0}^{J(K)-1}\sum_{ k_1,k_2 \in\Z^d}\sum_{e\in\bar
E} \beta_{j;k_1,k_2}^{(c,e)}\Psi_{j;k_1,k_2}^{(e)}(x,y)
\]
one can choose $(B=B(K,L))_{L\in\N}$ such that $\max_{x,y\in
[-B,B]^d}|\widetilde h^{(K)}_c(x,y)-\widetilde h^{(K,L)}_c(x,y)|=0$ and
$B(K,L)\longrightarrow_{L\to\infty} \infty$. This setting allows for
the approximations
\begin{eqnarray*}
E_1&\leq& C \tau_r^\delta\bigl[\E\bigl|H^{(L)}(X_k,X_l)\bigr|^{1/(1-\delta)}\I
_{(X_k^\prime,X_l^\prime)^\prime\in[-B,B]^d}\bigr]^{1-\delta}\leq C
\tau_r^\delta\bigl[P(X_1\notin[-B,B]^d)\bigr]^{1-\delta},\\
E_2&\leq& C \tau_r^\delta\bigl[P(X_1\notin[-B,B]^{2d})\bigr]^{1-\delta}.
\end{eqnarray*}
Analogously, it can be shown that
$\sup_{k\in\N}\E[H^{(L)}(X_1,X_{1+k})]^2\leq C P(X_1\notin[-B,B]^d)$.
Finally, we obtain
\[
\sup_{n\in\N}n^2 \E\bigl(U_{n,c}^{(K)}-U_{n,c}^{(K,L)}\bigr)^2
\leq C  \bigl[P(X_1\notin[-B,B]^d)\bigr]^{1-\delta}\Biggl[\sup_{n\in\N}\sum
_{r=1}^{n-1} (r+1)\tau_r^\delta\Biggr]\mathop{\longrightarrow}
_{L\to\infty} 0.
\]
Hence, the relations~(\ref{eq.l5.1}) and~(\ref{eq.l5.2}) hold.
\end{pf*}

\begin{pf*}{Proof of Lemma~\textup{\protect\ref{l.4}}}
In order to prove the assertion, we follow the lines of the proofs of
Lemma~\ref{l.1}, Lemma~\ref{l.5}, and Lemma~\ref{l.2} and carry out
some modifications.

In a first step, we reduce the problem to statistics with bounded
kernels $h_c$ defined in the proof of Lemma~\ref{l.1}. To this end, we
use the modified approximation
\begin{eqnarray*}
\bigl|H^{(c)}(x,y)-H^{(c)}(\bar x,\bar y)\bigr|
&\leq&[2 f(x,\bar x, y,\bar y)+g(x,\bar x)+g(y,\bar y)]  [\|x-\bar x\|
_{l_1}+\|y-\bar y\|_{l_1}]\\
&=:&f_1(x,\bar x, y,\bar y) [\|x-\bar x\|_{l_1}+\|y-\bar y\|_{l_1}],
\end{eqnarray*}
where $g$ is given by $g(x,\bar x):=\int_{\R^d}f(x,\bar x,z,z)P_X(\mathrm{d}z)$.
Under (A4)(i) H\"older's inequality yields
\begin{eqnarray*}
&&\E\bigl|H^{(c)}(Y_{k_1},Y_{k_2})-H^{(c)}(Y_{k_3},Y_{k_4})\bigr|\\
&&\quad\leq\Biggl(\E[f_1(Y_{k_1},Y_{k_2},Y_{k_3},Y_{k_4})]^{1/(1-\delta)}\sum
_{i=1}^4\|Y_{k_i}\|_{l_1}\Biggr)^{1-\delta} (\E\|Y_{k_1}-Y_{k_3}\|
_{l_1}+\E\|Y_{k_2}-Y_{k_4}\|_{l_1})^\delta
\end{eqnarray*}
for $Y_{k_i}  (k_i=1,\ldots,5, i=1,\ldots,4)$, as defined in (A4).
Plugging in this inequality into the calculations of the proof of
Lemma~\ref{l.1} yields $\sup_{n\in\N} n^2 \E(U_n-U_n^{(c)})^2
\longrightarrow_{c\to\infty}0.$

The next step contains the wavelet approximation of the bounded kernel
$h_c$. Defining~$h_c^{(K)}$ and $U_{n,c}^{(K)}$ as in the proof of
Lemma~\ref{l.5}, analogous to the proof of Lemma~\ref{l.2} there exists
a $C>0$ such that
%
\begin{eqnarray}\label{eq.t2.1}
&&\bigl|\widetilde h^{(K)}_c(\bar x,\bar y)-\widetilde h^{(K)}_c(x,y)\bigr|\nonumber\\
&&\quad\leq  f_1(x,\bar x,y, \bar y)[\|x-\bar x\|_{l_1}+\|y-\bar y\|
_{l_1}]+|H_2(\bar x,\bar y)-H_2(x,y)|
\nonumber
\\[-8pt]
\\[-8pt]
\nonumber
&&\quad\leq C f_1(x,\bar x,y, \bar y)[\|x-\bar x\|_{l_1}+\|y-\bar y\|
_{l_1}]\\
&&\qquad{} +\sum_{k_1,k_2\in\Z^d} \bigl(|\kappa_{k_1,k_2}(\bar x, \bar
y)|\bigl|\Phi_{J(K),k_1}(x)\Phi_{J(K),k_2}(y)-\Phi_{J(K),k_1}(\bar x)\Phi
_{J(K),k_2}(\bar y)\bigr|\bigr),\nonumber
\end{eqnarray}
where
$\kappa_{k_1,k_2}$ is given by
\[
\kappa_{k_1,k_2}(x, y):=
\int_{\R^d}\int_{\R^d}\Phi_{J(K),k_1}(u)\Phi
_{J(K),k_2}(v)[h_c(u,v)-h_c(x,y)]\, \mathrm{d}u \,\mathrm{d}v
\]
and $H_2$ is defined as in the proof of Lemma~\ref{l.2}.
In order to approximate the last summand of~(\ref{eq.t2.1}), we
distinguish again between the cases whether or not $(\bar x^\prime,
\bar y^\prime)^\prime\in\operatorname{supp }(\Phi_{J(K),k_1}\times \Phi_{J(K),k_2})$.
In the first case, an upper bound of order
\[
\mathrm{O}\Bigl(\max_{a_1,a_2\in[-S_\phi/2^{J(K)},S_\phi/2^{J(K)}]^d} f_1(\bar x,
\bar x+a_1,\bar y,\bar y+a_2)\Bigr)(\|\bar x-x\|_{l_1}+ \|\bar y-y\|_{l_1})
\]
can be obtained since
\begin{eqnarray*}
|\kappa_{k_1,k_2}(\bar x,\bar y)|&\leq&\frac{ S_\phi}{2^{J(K)}} \max
_{a_1,a_2\in[-S_\phi/2^{J(K)},S_\phi/2^{J(K)}]^d}f_1(\bar x, \bar
x+a_1,\bar y,\bar y+a_2)\\
&&{} \times\iint_{\R^d\times\R^d}\bigl|\Phi_{J(K),k_1}(u)\Phi
_{J(K),k_2}(v)\bigr| \,\mathrm{d}u \,\mathrm{d}v\\
&\leq& \mathrm{O}\bigl(2^{-J(K)(d+1)}\bigr)\max_{a_1,a_2\in[-S_\phi
/2^{J(K)},S_\phi/2^{J(K)}]^d}f_1(\bar x, \bar x+a_1,\bar y,\bar y+a_2).
\end{eqnarray*}
Here, $S_\phi$ denotes the length of the support of $\phi$. In the
second case, a decomposition similar to~(\ref{eq.l2.1}) can be employed
which leads to
the upper bound
\[
\mathrm{O}\Bigl(f_1(x,\bar x,y,\bar y)+\max_{a_1,a_2\in[-S_\phi/2^{J(K)},S_\phi
/2^{J(K)}]^d} f_1( x, x+a_1, y, y+a_2)\Bigr)(\|\bar x-x\|_{l_1}+ \|\bar y-y\|_{l_1}).
\]
Consequently, we get
\begin{eqnarray*}
\bigl|\widetilde h^{(K)}_c(\bar x,\bar y)-\widetilde h^{(K)}_c(x,y)\bigr|&\leq&
\mathrm{O}\Bigl( \max_{a_1,a_2\in[-S_\phi/2^{J(K)},S_\phi/2^{J(K)}]^d} f_1( x,
x+a_1, y, y+a_2)\\
&&\hphantom{\mathrm{O}\Bigl(}{} + \max_{a_1,a_2\in[-S_\phi/2^{J(K)},S_\phi/2^{J(K)}]^d} f_1(\bar x,
\bar x+a_1,\bar y,\bar y+a_2) \\
&&\hphantom{\mathrm{O}\Bigl(}{}+f_1(x,\bar x,y,\bar y)\Bigr)\times(\|\bar x-x\|_{l_1}+ \|\bar y-y\|
_{l_1})\\
&=:& f_2(x,\bar x, y, \bar y)(\|\bar x-x\|_{l_1}+ \|\bar y-y\|_{l_1}).
\end{eqnarray*}
This yields $|H^{(K)}(x,y)-H^{(K)}(\bar x, \bar y)|\leq f_3(x,\bar x,y
,\bar y)(\|x-\bar x\|_{l_1}+\|y-\bar y\|_{l_1})$ with $f_3(x,\bar x,y,\allowbreak
\bar y)=2 f_2(x,\bar x,y ,\bar y)+\int_{\R^d} f_2(x,\bar
x,z,z)P_X(\mathrm{d}z)+\int_{\R^d} f_2(z,z,\bar y,y)P_X(\mathrm{d}z)$.
Note that under~(A4)(i), $\E[ f_3(Y_i,Y_j,Y_k,Y_l)]^\eta(\|Y_i\|
_{l_1}+\|Y_j\|_{l_1}+\|Y_k\|_{l_1}+\|Y_l\|_{l_1})<\infty$ if $J(K)$ is
sufficiently large. Thus, we have
\[
\E\bigl|H^{(K)}(Y_{k_1},Y_{k_2})-H^{(K)}(Y_{k_3},Y_{k_4})\bigr|
\leq C (\E\|Y_{k_1}-Y_{k_3}\|_{l_1}+\E\|Y_{k_2}-Y_{k_4}\|_{l_1})^\delta
\]
for $Y_{k_i}  (k_i=1,\ldots,5, i=1,\ldots,4)$, as defined in (A4).
Moreover, Lemma~\ref{l.3} remains valid with $g=h_c$.
Therefore, one can follow the lines of the proof of Lemma~\ref{l.4} and
plug in the inequality above. This procedure leads to $\sup_{n\in\N
}n^2 \E(U_{n,c}-U_{n,c}^{(K)})^2 \longrightarrow_{K\to\infty} 0$.

In the third step of the proof, we verify $\sup_{n\in\N} n^2 \E
(U_{n,c}^{(K)}-U_{n,c}^{(K,L)})^2 \longrightarrow_{L\to\infty}0$.
For this purpose, it suffices to plug in a modified approximation of
$H^{(L)}(x,y)-H^{(L)}(\bar x,\bar y)$ into the second part of the proof
of Lemma~\ref{l.5}.
Lipschitz continuity of $h_c^{(K,L)}$ implies
\[
\bigl|H^{(L)}(x,y)-H^{(L)}(\bar x,\bar y)\bigr|\leq f_4(x,\bar x,y ,\bar y)[\|
x-\bar x\|_{l_1}+\|y-\bar y\|_{l_1}]\vspace*{-2pt}
\]
with $f_4(x,\bar x,y ,\bar y)=C+f_3(x,\bar x,y ,\bar y)$.
Since, $f_4$ satisfies the moment assumption of (A4)(i) with $A=0$ for
sufficiently large $J(K)$, we obtain
\[
\E\bigl|H^{(L)}(Y_{k_1},Y_{k_2})-H^{(L)}(Y_{k_3},Y_{k_4})\bigr|\leq C [\E(\|
Y_{k_1}-Y_{k_3}\|_{l_1}+\|Y_{k_2}-Y_{k_4}\|_{l_1})]^\delta.\vspace*{-2pt}
\]
Hence, $\sup_{n\in\N} n^2 \E(U_{n,c}^{(K)}-U_{n,c}^{(K,L)})^2
\longrightarrow_{L\to\infty}0$.
Summing up the three steps yields
\[
\lim_{c\to\infty}\limsup_{K\to\infty}\limsup_{L\to\infty}\sup_{n\in\N
}n^2 \E\bigl(U_n-U_{n,c}^{(K,L)}\bigr)^2=0.
\]
\upqed\vspace*{-2pt}\end{pf*}

\begin{pf*}{Proof of Lemma~\textup{\protect\ref{l.6}}}
A positive variance of $Z$ implies the existence of constants $V>0$ and
$c_0>0$ such that for every $c\geq c_0$ we can find a $K_0\in\N$ such
that for every $K\geq K_0$ there is an $L_0$ with $\operatorname{var}(Z^{(K,L)}_c)\geq
V, \forall L\geq L_0.$ Moreover, uniform equicontinuity of the
distribution functions of $(((Z^{(K,L)}_c)_{L})_{K})_c$ yields the
desired property of $Z$. By matrices-based notation of $Z_c^{(K,L)}$,
we obtain
\[
Z^{(K,L)}_c =C^{(K,L)}+\sum_{k_1,k_2=1}^{M(K,L)}\gamma
_{k_1,k_2}^{(c,K,L)}Z_{k_1}^{(K,L)} Z_{k_2}^{(K,L)}=C^{(K,L)}+\bigl[\bar
Z^{(K,L)}\bigr]^\prime\Gamma^{(K,L)}_c\bar Z^{(K,L)},\vspace*{-2pt}
\]
with a constant $C^{(K,L)}$, a symmetric matrix of coefficients $\Gamma
^{(K,L)}_c$, and a normal vector
$\bar Z^{(K,L)}=(Z_1^{(K,L)},\ldots, Z_{M(K,L)}^{(K,L)})^\prime$.
Hence, $Z^{(K,L)}_c-C^{(K,L)}$ can be rewritten as follows:
\begin{eqnarray*}
Z^{(K,L)}_c-C^{(K,L)}&\stackrel{d}{=}& \bar Y^\prime
\bigl[U_c^{(K,L)}\bigr]^\prime\Lambda^{(K,L)}_c U_c^{(K,L)}\bar Y
=Y^\prime \Lambda^{(K,L)}_c Y\\[-2pt]
&=&\sum_{k=1}^{M(K,L)}\lambda_k^{(c,K,L)}Y_k^2.\vspace*{-2pt}
\end{eqnarray*}
Here $U_c^{(K,L)}$ is a certain orthogonal matrix, $\Lambda
^{(K,L)}_c:=\operatorname{diag}(\lambda_1^{(c,K,L)}, \ldots,\lambda
_{M{(K,L)}}^{(c,K,L)})$ with $|\lambda_1^{(c,K,L)}|\geq\cdots\geq|\lambda
_{M{(K,L)}}^{(c,K,L)}|$, and $\bar Y$ as well as $ Y$ are multivariate
standard normally distributed random vectors. For notational
simplicity, we suppress the upper index $(c,K,L)$ in the sequel.
Due to the above choice of the triple~$(c,K,L)$, either $\sum
_{k=1}^{4}(\lambda_k)^2$ or $\sum_{k=5}^{M(K,L)}(\lambda
_k)^2$ is bounded from below by $V/4$.
In the first case, $\lambda_1\geq\sqrt{V/16}$ holds true which implies
\[
P\bigl(Z^{(K,L)}_c\in[x-\varepsilon,x+\varepsilon]\bigr)\leq\int
_{0}^{2\varepsilon} f_{\lambda_1Y_1^2}(t)\, \mathrm{d}t
\leq P(Y_1^2\leq2\varepsilon)\max\biggl\{1,\frac{4}{\sqrt V}\biggr\}\qquad
 \forall  x\in\R.\vspace*{-2pt}
\]
Here, the first inequality results from the fact that convolution
preserves the continuity properties of the smoother function.\vadjust{\goodbreak}
In the opposite case, that is, $\sum_{k=5}^{M(K,L)}(\lambda_k
)^2\geq V/4$, it is possible to bound the uniform norm of the density
function of $Z^{(K,L)}_c$ by means of its variance. To this end, we
first consider the characteristic function $\varphi_{Z^{(K,L)}_c}$ of
$Z^{(K,L)}_c$ and assume w.l.o.g.~that $M(K,L)$ is divisible by 4.
Defining a sequence $(\mu_k)_{k=1}^{M(K,L)/4}$ by $\mu_k=\lambda_{4k}$
for $k\in\{1,\ldots,M(K,L)/4\}$ allows for the approximation:
\begin{eqnarray*}
\bigl| \varphi_{Z^{(K,L)}_c}(t) \bigr|
& =& \Biggl\{ \prod_{j=1}^{M(K,L)} ( 1  +  [2\lambda_j t]^2
) \Biggr\}^{-1/4} \leq\Biggl\{ \prod_{j=1}^{M(K,L)/4} ( 1  +
[2\mu_j t]^2 ) \Biggr\}^{-1} \\[-3pt]
& \leq &\frac{1}{1  +  4(\mu_1^2+\cdots+\mu_{M(K,L)/4}^2)  t^2}.
\end{eqnarray*}
By inverse Fourier transform, we obtain the following result concerning
the density function of $Z^{(K,L)}_c$:
\begin{eqnarray*}
\bigl\| f_{Z^{(K,L)}_c} \bigr\|_\infty
& \leq& \frac{1}{2\uppi} \| \varphi_{Z^{(K,L)}_c} \|_1 \leq \frac{1}{2\uppi
} \int_{-\infty}^\infty
\frac{1}{1 +  (2\sqrt{\mu_1^2+\cdots+\mu_{M(K,L)/4}^2} t)^2}  \, \mathrm{d}t \\[-3pt]
& =& \frac{1}{\sqrt{\mu_1^2+\cdots+\mu_{M(K,L)/4}^2}} \frac{1}{2\uppi}
\int_0^\infty\frac{1}{1+u^2}   \mathrm{d}u \\[-3pt]
&\leq&\frac{1}{2  \sqrt{4(\mu_1^2+\cdots+\mu_{M(K,L)/4-1}^2})}\\[-3pt]
&\leq&\frac{1}{2  \sqrt{\lambda_5^2+\cdots+\lambda_{M(K,L)}^2}}
\leq\frac{1}{\sqrt{V}}.
\end{eqnarray*}
Thus, $P(Z^{(K,L)}_c\in[x-\varepsilon,x+\varepsilon])\leq2\varepsilon/\sqrt
{V}$ which completes the studies of the case $\sum_{k=5}^{M(K,L)}
(\lambda_k)^2>V/4$ and finally yields the assertion.\vspace*{-2pt}
\end{pf*}

\begin{pf*}{Proof of Lemma~\textup{\protect\ref{l.7}}}
This result is an immediate consequence of Theorem~\ref{t.3}.\vspace*{-2pt}
\end{pf*}

\section*{Acknowledgements}\vspace*{-2pt}
The author is grateful to Michael H.~Neumann for his constructive
advice and fruitful discussions. She also thanks an anonymous referee
for helpful comments that led to an improvement of the paper.
This research was funded by the German Research Foundation DFG (project:
NE~\mbox{606/2-1}).\vspace*{-2pt}

%

\printhistory

\end{document}